\documentclass{amsproc}
\usepackage{amssymb}
\usepackage{graphicx}
\usepackage{amsmath}
\usepackage{url,paralist}
\usepackage{anysize}

\newtheorem{theorem}{Theorem}[section]
\newtheorem{lemma}[theorem]{Lemma}
\theoremstyle{definition}
\newtheorem{definition}[theorem]{Definition}
\newtheorem{proposition}[theorem]{Proposition}
\newtheorem{problem}[theorem]{Problem}

\newtheorem{conjecture}[theorem]{Conjecture}
\newtheorem{corollary}[theorem]{Corollary}
\newtheorem{notation}[theorem]{Notation}
\theoremstyle{remark}
\newtheorem{remark}[theorem]{Remark}
\numberwithin{equation}{section}

\begin{document}
\title{Using Equivariant Obstruction Theory in Combinatorial Geometry}
\author{Pavle V. M. Blagojevi\'{c}}
\address{Mathematical Institute SANU, Belgrade, Serbia}
\email{pavleb@mi.sanu.ac.yu}
\thanks{Supported by the grant 144018 of the Serbian Ministry of Science and
Environment.}
\author{Aleksandra S. Dimitrijevi\'{c} Blagojevi\'{c}}
\address{Faculty for Agriculture, Zemun, Serbia}
\email{vxdig@beotel.net}
\thanks{Supported by the grant 144018 of the Serbian Ministry of Science and
Environment. \\
Our deepest thanks go to referees, referee of the papers \cite{Bl-Vr-Ziv},
\cite{Bl}, professors Paul Goerss, John McCleary, Carsten Schultz and G\"{u}%
nter Ziegler. They provided valuable assistance, ideas and motivation.}
\subjclass{Primary 55S35, 52C35; Secondary 68U05}
\keywords{Equivariant Obstruction Theory, Mass partition, Graph Coloring}

\begin{abstract}
\noindent A significant group of problems coming from the realm of
Combinatorial Geometry can be approached fruitfully through the use of
Algebraic Topology. From the first such application to Kneser's problem in
1978 by Lov\'{a}sz \cite{Lovasz} through the solution of the Lov\'{a}sz
conjecture \cite{Babson-Kozlov}, \cite{Carsten}, many methods from Algebraic
Topology have been developed. Specifically, it appears that the
understanding of equivariant theories is of the most importance. The
solution of many problems depends on the existence of an elegantly
constructed equivariant map. For example, the following problems were
approached by discussing the existence of appropriate equivariant maps. A
variety of results from algebraic topology were applied in solving these
problems. The methods used ranged from well known theorems like Borsuk-Ulam
and Dold theorem to the integer / ideal-valued index theories. Recently
equivariant obstruction theory has provided answers where the previous
methods failed. For example, in papers \cite{Z4} and \cite{Bl} obstruction
theory was used to prove the existence of different mass partitions. In this
paper we extract the essence of the equivariant obstruction theory in order
to obtain an effective \emph{general position map} scheme for analyzing the
problem of existence of equivariant maps. The fact that this scheme is
useful is demonstrated in this paper with three applications:

\textbf{(A)} a "half-page" proof of the Lov\'{a}sz conjecture due to Babson
and Kozlov \cite{Babson-Kozlov} (one of two key ingredients is Carsten's map
\cite{Carsten}),

\textbf{(B)} a generalization of the result of V. Makeev \cite{Mak} about
the sphere $S^{2}$ measure partition by $3$-planes (Section \ref{Result1}),
and

\textbf{(C)} the new $(a,b,a)$, class of $3$-fan $2$-measures partitions
(Section \ref{Result2}).

\noindent These three results, sorted by complexity, share the spirit of
analyzing equivariant maps from spheres to complements of arrangements of
subspaces.
\end{abstract}

\maketitle






\section{Equivariant obstruction theory}

\noindent The basic concept of any obstruction theory is to produce an
invariant associated to a specific construction in such a way that the
nature of the invariant determines whether the construction can or can not
be performed. An (equivariant) obstruction theory considers two basic
problems. For a finite group $G$, consider a relative $G$-cellular complex $%
(X,A)$ such that the $G$-action on $X\backslash A$ is \textbf{free.} Let $Y$
be a $G$-space.

\noindent \textit{Extension problem. }Let $f:A\rightarrow Y$ be a $G$-map.
Is there a $G$-map $F:X\rightarrow Y$ such that $f=F\circ i$? Here $%
i:A\rightarrow X$ denotes the inclusion.

\noindent \textit{Homotopy problem.} Let $f_{0}:X\rightarrow Y$ and $%
f_{1}:X\rightarrow Y$ be $G$-maps such that there is a $G$-homotopy $%
h:I\times A\rightarrow Y$ from $f_{0}|_{A}$ to $f_{1}|_{A}$. Is there a $G$%
-homotopy $H:I\times X\rightarrow Y$ which extends $h$, i.e, $%
H|_{\{0\}\times X}=f_{0}$, $H|_{\{1\}\times X}=f_{1}$ and $H|_{I\times A}=h$?

\noindent The answer which obstruction theory provides is a sequence of
obstruction elements living in equivariant cohomology. For the details about
(equivariant) obstruction theory one can consult the expositions in \cite[%
Section II. 3]{Dieck87}, \cite[Chapter 7]{Davis} and \cite[Section V. 5]{Wh}.

\subsection{Equivariant homology and cohomology.}

\noindent Let $(X,A)$ be a relative $G$-cellular complex with a free action
on $X\backslash A$. Let $C_{\ast }(X,A)$ denote the integral cellular chain
complex. The cellular free $G$-action on every skeleton of $X\backslash A$
induces a free $G$-action on the chain complex $C_{\ast }(X,A)$. Therefore,
the chain complex $C_{\ast }(X,A)$ is actually a chain complex of free $%
\mathbb{%
\mathbb{Z}
}[G]$-modules.

\begin{definition}
Let $(X,A)$ be a relative $G$-cellular complex such that the $G$-action on $%
X\backslash A$ is free. Let $M$ be a $\mathbb{%
\mathbb{Z}
}[G]$-module.

(A) The chain complex%
\begin{equation}
\mathfrak{C}_{\ast }^{G}(X,A;M)=C_{\ast }(X,A)\otimes _{\mathbb{%
\mathbb{Z}
}[G]}M  \label{eq:0}
\end{equation}%
is the equivariant chain complex of $(X,A)$ with coefficients in $M$, and
its homology $\mathfrak{H}_{\ast }^{G}(X,A;M)$ is called the \textbf{%
equivariant homology} of $(X,A)$ with coefficients in $M$.

(B) The cochain complex
\begin{equation}
\mathfrak{C}_{G}^{\ast }(X,A;M)=\mathrm{Hom}_{\mathbb{%
\mathbb{Z}
}[G]}(C_{\ast }(X,A),M)  \label{eq:1}
\end{equation}%
is the equivariant cochain complex of $(X,A)$ with coefficients in $M$, and
its homology $\mathfrak{H}_{G}^{\ast }(X,A;M)$ is called the \textbf{%
equivariant cohomology} of $(X,A)$ with coefficients in $M$.
\end{definition}

\noindent A central example of a $\mathbb{%
\mathbb{Z}
}[G]$-module in our applications will come from a path connected, $n$-simple
$G$-space $Y$. Being $n$\textbf{-simple} means that $\pi _{1}(Y,y_{0})$ acts
trivially on $\pi _{n}(Y,y_{0})$ for every $y_{0}\in Y$. The action of $G$
on $Y$ easily produces an action on the set of homotopy classes $[S^{n},Y]$.
Since $Y$ is $n$-simple, the action can be extended to the homotopy group $%
\pi _{n}Y\cong \lbrack S^{n},Y]$. Consequently there exists a $\mathbb{%
\mathbb{Z}
}[G]$-module structure on $\pi _{n}Y$.

\noindent The equivariant homology and cohomology groups of $(X,A)$ can be
interpreted as the ordinary homology and cohomology groups of the quotient
pair $(X/G,A/G)$ with appropriate local coefficients \cite[p.112]{Dieck87}.
For $M$ interpreted as a local coefficient system $\tilde{M}$ on $%
(X/G)\backslash (A/G)$:%
\begin{equation}
\mathfrak{H}_{n}^{G}(X,A;M)\cong H_{n}(X/G,A/G;\tilde{M}),  \label{eq:2}
\end{equation}%
\begin{equation}
\mathfrak{H}_{G}^{n}(X,A;M)\cong H^{n}(X/G,A/G;\tilde{M}).  \label{eq:4}
\end{equation}

\subsection{The exact obstruction sequence.}

\noindent Let $n\geq 1$ be a fixed integer and $Y$ a path-connected $n$%
-simple $G$-space. For every $G$-relative cell complex $(X,A)$ with \textbf{%
free} action of $G$ on $X\backslash A$, there exists an obstruction exact
sequence%
\begin{equation}
\lbrack X_{n+1},Y]_{G}\longrightarrow \mathrm{im}\left(
[X_{n},Y]_{G}\rightarrow \lbrack X_{n-1},Y]_{G}\right) \overset{\mathfrak{o}%
_{G}^{n+1}}{\longrightarrow }\mathfrak{H}_{G}^{n+1}(X,A;\pi _{n}Y),
\label{ObstructionExactSequence}
\end{equation}%
where $X_{k}$ denotes the $k$-skeleton of $X$. The sequence is natural in $X$
and $Y$. This exact sequence should be understood in the following way:

(A) Every $G$-map on the $(n-1)$-skeleton $f:X_{n-1}\rightarrow Y$ which can
be equivariantly extended to the $n$-skeleton $f:X_{n}\rightarrow Y$ defines
an unique element $\mathfrak{o}_{G}^{n+1}(f)$ living in $\mathfrak{H}%
_{G}^{n+1}(X,A;\pi _{n}Y)$, called the \textbf{obstruction element};

(B) The exactness of the sequence means that the obstruction element $%
\mathfrak{o}_{G}^{n+1}(f)$ is zero if and only if there is a map in the
homotopy class of the restriction $f|_{X_{n-1}}$ which can be extended to
the $(n+1)$-skeleton $X_{n+1}$.

\noindent The obstruction element can be introduced on the cochain level by
a universal geometrical construction. Let $[h]\in \lbrack X_{n},Y]_{G}$, let
$\varphi :(D^{n+1},S^{n})\rightarrow (X_{n+1},X_{n})$ be an attaching map
and $e\in C_{n+1}(X,A)$ the generator associated to the cell $\varphi $.
Then the \textbf{obstruction cochain}\textit{\ }$\mathfrak{o}%
_{G}^{n+1}(h)\in \mathfrak{C}_{G}^{n+1}(X,A;\pi _{n}Y)$ of the map $h$ is
defined on $e$ (with a little abuse of the notation) by%
\begin{equation}
\mathfrak{o}_{G}^{n+1}(h)(e)=[h\circ \varphi ]\in \lbrack S^{n},Y].
\label{DefObstructionCocycle}
\end{equation}%
It can be proved that the cohomology class of the obstruction cocycle is the
obstruction element defined by the exact sequence (\ref%
{ObstructionExactSequence}).

\begin{figure}[tbh]
\centering\includegraphics[scale=0.70]{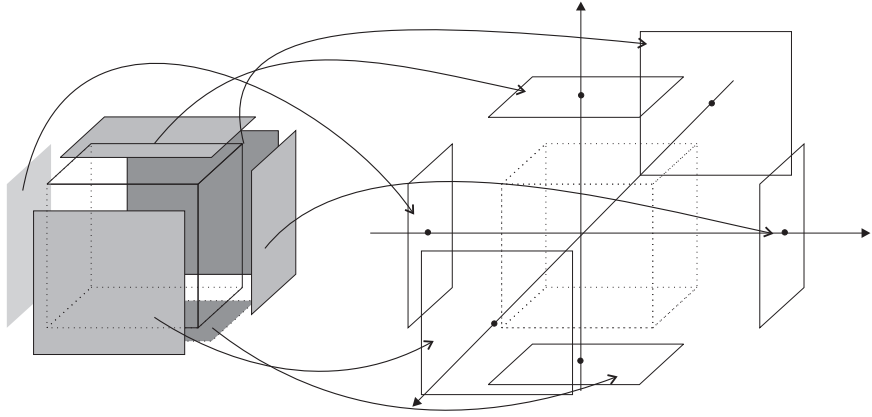}
\caption{{\protect\small The obstruction cocycle for the map of a cube in
the complement of the coordinate lines arrangement.}}
\label{fig:Fig1}
\end{figure}

\subsection{The primary obstruction.}

\noindent Let $n\geq 1$ be a fixed integer. The following proposition holds
for an $(n-1)$-connected, $n$-simple $G$-space $Y$.

\begin{proposition}
\label{prop:PrimalObstruction}Let $(X,A)$ be a relative $G$-cell complex
with the action of $G$ free on $X\backslash A$, and $f:A\rightarrow Y$ any $%
G $-map.

(A) There exists a $G$-map $h:X_{n}\rightarrow Y$ extending $f$, i.e., $%
h|_{A}=f$.

(B) Every two $G$-extensions $h$ and $k:X_{n}\rightarrow Y$ of $f$ are $G$%
-homotopic $\mathrm{rel~}A$ on $X_{n-1}$, that is
\begin{equation*}
\mathrm{im}\left( [X_{n},Y]_{G}\rightarrow \lbrack X_{n-1},Y]_{G}\right)
=\{\ast \}.
\end{equation*}

(C) If $H:I\times A\rightarrow Y$ is a $G$-homotopy between $f:A\rightarrow
Y $ and $g:A\rightarrow Y$, and $h,~k:X_{n}\rightarrow Y$ are the extensions
of $f$ and $g$, then there is a $G$-homotopy $K:I\times X_{n-1}\rightarrow Y$
extending $H$ between $h|_{X_{n-1}}$ and $k|_{X_{n-1}}$.
\end{proposition}

\begin{proof}
(A) We extend $f=f_{0}$ starting with the $0$-skeleton and go up to the $n$%
-skeleton. The obstruction for lifting $f_{r}:X_{r}\rightarrow Y$ from the $%
r $-skeleton to the $(r+1)$-skeleton lies in $\mathfrak{C}_{G}^{r+1}(X,\pi
_{r}Y)$. Since $Y$ is $(n-1)$-connected and $n$-simple,
\begin{equation*}
\pi _{r}Y=0\text{ for all }1\leq r\leq n-1~\Rightarrow ~\mathfrak{C}%
_{G}^{r+1}(X,\pi _{r}Y)=0\text{ for all }1\leq r\leq n-1\text{.}
\end{equation*}%
Hence, there is a $G$-map $h:X_{n}\rightarrow Y$ extending $f$.

\noindent (B) Let $h,~k:X_{n}\rightarrow Y$ be $G$-extensions of $f$. A $G$%
-map%
\begin{equation*}
K:\left( \{0\}\times X_{n}\right) \cup \left( \{1\}\times X_{n}\right) \cup
\left( I\times A\right) \rightarrow Y
\end{equation*}%
can be defined by
\begin{equation*}
\begin{array}{ll}
K|_{\{0\}\times X_{n}}=h, & K|_{\{1\}\times X_{n}}=k, \\
K|_{\{t\}\times A}=f\text{,} & \text{for every }t\in I.%
\end{array}%
\end{equation*}%
Now consider a relative $G$-cell complex $(I\times X_{n},\left( \{0\}\times
X_{n}\right) \cup \left( \{1\}\times X_{n}\right) \cup \left( I\times
A\right) )$ and extend $K$. The assumptions on $Y$ provide a $G$-map $H$
from the $n$-skeleton $\left( \{0\}\times X_{n}\right) \cup \left(
\{1\}\times X_{n}\right) \cup \left( I\times A\right) $ of $I\times X_{n}$
to $Y$ which extends $K$. The restriction $H|_{I\times X_{n-1}}$ is the
required $G$-homotopy $\mathrm{rel~}A$.

\noindent (C) Considering a relative $G$-cell complex $(I\times
X_{n-1},\left( \{0\}\times X_{n}\right) \cup \left( \{1\}\times X_{n}\right)
\cup \left( I\times A\right) )$ instead of $(X,A)$ the statement becomes a
direct consequence of the property (A).
\end{proof}

\noindent When%
\begin{equation*}
\mathrm{im}\left( [X_{n},Y]_{G}\rightarrow \lbrack X_{n-1},Y]_{G}\right)
=\{\ast \},
\end{equation*}%
the obstruction sequence (\ref{ObstructionExactSequence}) becomes

\begin{equation}
\lbrack X_{n+1},Y]_{G}\longrightarrow \{\ast \}\overset{\mathfrak{o}%
_{G}^{n+1}}{\longrightarrow }\mathfrak{H}_{G}^{n+1}(X,\pi _{n}Y).
\label{PrimaryObstrSequence}
\end{equation}%
The element $\mathfrak{o}_{G}^{n+1}(\ast )\in \mathfrak{H}_{G}^{n+1}(X,\pi
_{n}Y)$ is called the \textbf{primary obstruction} and does not depend on
the map of the $n$-th skeleton.

\begin{corollary}
If $\circ $ and $\ast $ are $G$-actions on $S^{n}$ and $\circ $ is free,
then there exists a $G$-map $f:S^{n}\rightarrow S^{n}$ such that
\begin{equation*}
f(g\circ x)=g\ast f(x)
\end{equation*}%
for all $g\in G$ and $x\in S^{h}$.
\end{corollary}

\noindent The Corollary is a direct consequence of the previous proposition
statement (A).

\subsection{Equivariant Poincar\'{e} duality.}

\noindent Let $X$ be a compact $n$-dimensional free $G$-manifold. Then there
is a version of the Poincar\'{e} duality isomorphism for the equivariant
homology and cohomology of$\mathfrak{\ }X$ with the coefficients in any $G$%
-module $M$.

\begin{theorem}
\label{Th:EqPoincare}Let $X$ be a compact $n$-dimensional, simply connected,
free $G$-manifold, $\mathcal{Z}$, a $G$-module, $H_{n+1}(X,%
\mathbb{Z}
)\cong
\mathbb{Z}
$, and $M$ a $G$-module. Then there is an isomorphism of the groups%
\begin{equation*}
\mathfrak{H}_{G}^{k}(X;M)\cong \mathfrak{H}_{n-k}^{G}(X;M\otimes \mathcal{Z})
\end{equation*}%
for every $k\in \{0,...,n\}$.
\end{theorem}

\begin{proof}
When $M$ is interpreted as the local coefficient system $\tilde{M}$, the
isomorphism (\ref{eq:4}) implies:
\begin{equation*}
\mathfrak{H}_{G}^{k}(X;M)\cong H^{k}(X/G;\tilde{M}).
\end{equation*}%
From \cite[Theorem 4.18. p. 196]{Kaw} and \cite[Proposition 1.40. p. 70]%
{hatcher} we know that $X/G$ is a compact manifold with fundamental group $%
\pi _{1}(X/G)\cong G$. Since all the assumptions for applying Theorem 2.1
\cite[p.23]{Wall} are met,
\begin{equation*}
H^{k}(X/G;\tilde{M})\cong H_{n-k}^{t}(X/G;\tilde{M}).
\end{equation*}
Careful reading of the definition of the modified homology $H_{\ast }^{t}($.$%
)$, \cite[p. 21]{Wall} implies that%
\begin{equation*}
H_{n-k}^{t}(X/G;\tilde{M})=H_{n-k}(X/G;\tilde{M}\otimes \mathcal{Z}).
\end{equation*}%
The isomorphism (\ref{eq:2}) concludes the proof%
\begin{equation*}
H_{n-k}(X/G;\tilde{M}\otimes \mathcal{Z})\cong \mathfrak{H}%
_{n-k}^{G}(X;M\otimes \mathcal{Z}).
\end{equation*}
\end{proof}

\subsection{The existence of a $G$-map $M\rightarrow W\backslash \Sigma $
from a manifold to a complement.}

\noindent Let $M$ be a connected, $(n+1)$-dimensional, compact free $G$%
-manifold, $W$ a $d$-dimensional smooth $G$-manifold, and $\Sigma $ the
union of a finite $G$-invariant arrangement $\mathcal{S}=\{S_{i}|i\in I\}$
of the $(d-n-1)$-dimensional smooth submanifolds. Let us also assume that

(A) the complement $W\backslash \Sigma $ is $n$-simple, paracompact space,

(B) the complement $W\backslash \Sigma $ is $(n-1)$-connected,

(C) the tangent spaces of the submanifolds $S_{i}$ in any mutual
intersection point do not coincide, and

(D) $H_{n}(W,%
\mathbb{Z}
)=0$.

\noindent \textit{The question we consider is whether there is a }$G$-%
\textit{map }$M\rightarrow W\backslash \Sigma $\textit{.}

\subparagraph{\textbf{1.}}

\noindent Since the complement $W\backslash \Sigma $ is $(n-1)$-connected by
assumption, the problem of the existence of a $G$-map $M\rightarrow
W\backslash \Sigma $ depends only on primary obstruction. The obstruction
exact sequence, like in (\ref{PrimaryObstrSequence}), has the form%
\begin{equation*}
\lbrack M,W\backslash \Sigma ]_{G}\longrightarrow \{\ast \}\overset{%
\mathfrak{o}_{G}^{n+1}}{\longrightarrow }\mathfrak{H}_{G}^{n+1}(M,\pi
_{n}\left( W\backslash \Sigma \right) ).
\end{equation*}%
The assumptions (A) and (B) on the space $W\backslash \Sigma $ and the
Hurewicz theorem imply that $\pi _{n}\left( W\backslash \Sigma \right) \cong
H_{n}(W\backslash \Sigma ,%
\mathbb{Z}
)$ as $G$-modules. Thus, the obstruction element $\mathfrak{o}%
_{G}^{n+1}(\ast )$ lives in the group $\mathfrak{H}_{G}^{n+1}(M,H_{n}(W%
\backslash \Sigma ,%
\mathbb{Z}
))$, where the natural $G$-structure on $H_{n}(W\backslash \Sigma ,%
\mathbb{Z}
)$ is assumed.

\noindent Consider the equivariant cohomology group $\mathfrak{H}%
_{G}^{n+1}(M,H_{n}(W\backslash \Sigma ,%
\mathbb{Z}
))$. Since $M$ is a free $G$-manifold by assumption, Theorem \ref%
{Th:EqPoincare} provides an isomorphism
\begin{equation*}
\mathfrak{H}_{G}^{n+1}(M,H_{n}(W\backslash \Sigma ,%
\mathbb{Z}
))\cong \mathfrak{H}_{0}^{G}(M,H_{n}(W\backslash \Sigma ,%
\mathbb{Z}
)\otimes \mathcal{Z})
\end{equation*}%
where $\mathcal{Z}$ is the $G$-module $H_{n+1}(M,%
\mathbb{Z}
)\cong
\mathbb{Z}
$. The isomorphism (\ref{eq:2}) implies that%
\begin{equation*}
\mathfrak{H}_{0}^{G}(M,H_{n}(W\backslash \Sigma ,%
\mathbb{Z}
)\otimes \mathcal{Z})\cong H_{0}\left( M/G;H_{n}(W\backslash \Sigma ,%
\mathbb{Z}
)\otimes \mathcal{Z}\right) .
\end{equation*}%
Since $G\cong \pi _{1}(M/G)$ an application of Proposition 5.14. of \cite[%
p.107]{Davis} (or alternatively \cite[Exercise 1, p.44]{Brown} and \cite[%
(1.5), p.57]{Brown}) provides an isomorphism
\begin{equation}
H_{0}\left( M/G;H_{n}(W\backslash \Sigma ,%
\mathbb{Z}
)\otimes \mathcal{Z}\right) \cong \left( H_{n}(W\backslash \Sigma ,%
\mathbb{Z}
)\otimes \mathcal{Z}\right) _{G}\text{.}  \label{CoinvariantIsomorphism}
\end{equation}%
Thus the obstruction element lives in a group of coinvariants of the first
non-trivial reduced homology group of the target space $W\backslash \Sigma $%
,
\begin{equation*}
\mathfrak{o}_{G}^{n+1}(\ast )\in \mathfrak{H}_{G}^{n+1}(M,H_{n}(W\backslash
\Sigma ,%
\mathbb{Z}
))\cong \left( H_{n}(W\backslash \Sigma ,%
\mathbb{Z}
)\otimes \mathcal{Z}\right) _{G}
\end{equation*}

\subparagraph{\textbf{2.}}

\noindent The situation where the primary obstruction is the only
obstruction, as in our case, has the advantage of not depending on the
particular $G$-map on the $n$-th skeleton.

\begin{definition}
\label{Def:MapInGeneralPosition}Let a $G$-map $f:M\rightarrow W$ satisfy the
following conditions

\noindent \textit{(A)} $f(M_{n})\subset W\backslash \Sigma $,

\noindent \textit{(B) }$f(M)\cap \Sigma $ is a finite set of points,

\noindent \textit{(C) }$(\forall x\in f(M)\cap \Sigma )(\forall S$ $\in
\mathcal{S})$ $\{x\}=S\cap f(M)~\Rightarrow ~$ intersection is transversal
at $x$,

\noindent \textit{(D) }$(\forall x\in f(M)\cap \Sigma )(\forall
S_{1},S_{2}\in \mathcal{S})~\{x\}=S_{1}\cap S_{2}\cap f(M)~\Rightarrow ~$%
\textrm{codim}$_{S_{1}}\left( S_{1}\cap S_{2}\right) =1$.

\noindent We then say that $f:M\rightarrow W$ is a map in \textbf{general
position} with respect to $\Sigma $.
\end{definition}

\noindent Condition (D) allows the intersection points $f(M)\cap \Sigma $ to
belong to the lower strata of the arrangement $\mathcal{S}$. This forces the
introduction of broken point classes along point classes as possible results
of evaluation of the obstruction cocycle.

\subparagraph{\textbf{3.}}

\noindent The notion of point and broken point classes was introduced in
\cite{Bl} and \cite{Bl-Vr-Ziv} for the complements of arrangements of linear
spaces. We extend this definition to the present setting. Consider $x\in
\Sigma $. There are elements $S_{1},..,S_{k}$ in $\mathcal{S}$ such that $%
x\in S_{1}\cap ..\cap S_{k}$ and \textrm{codim}$_{S_{i}}\left( S_{1}\cap
..\cap S_{k}\right) =1$. Let $D_{1},..,D_{k}$ denote disks in fibers at the
point $x$ of tubular neighborhoods of the submanifolds $S_{1},..,S_{k}$ such
that for all $i,j\in \{1,..,k\}$we have $S_{i}\cap D_{j}=\{x\}$. The
smoothness assumptions on $W$ and elements of the arrangement $\mathcal{S}$
guarantee the existence of the above construction \cite[Theorem 11.14. p.100]%
{Bredon}. The fundamental class of the pair $(D_{i},\partial D_{i})$
determines a homology class in $H_{n+1}(W,W\backslash \Sigma ;\mathbb{%
\mathbb{Z}
})$ that we denote by $[x,D_{i}]$ and call the\textit{\ }\textbf{point class}
of $x$ determined by $D_{i}$. The homotopy axiom implies that $[x,D_{i}]$
does not change if $x$ is moved inside the connected component of $%
S_{i}\setminus \bigcup \{S\neq S_{i}|S\in \mathcal{S\}}$. Because of the
assumption (D), at the beginning of the section, that $H_{n}(W,%
\mathbb{Z}
)=0$, the epimorphism $H_{n+1}(W,W\backslash \Sigma ;\mathbb{%
\mathbb{Z}
})\rightarrow H_{n}(W\backslash \Sigma )$ determines the class $\left\Vert
x,D_{i}\right\Vert :=\partial \lbrack x,D_{i}]$ which is also called the%
\textbf{\ point class}\textit{\ }of\emph{\ }$x$ determined by $D_{i}$. Thus
all point classes $\left\Vert x,D_{i}\right\Vert $ are born as $[x,D_{i}]$
classes, but $\left\Vert x,D_{i}\right\Vert $ can be zero while $%
[x,D_{i}]\neq 0$.

\begin{proposition}
\label{prop:SingularSet}Consider a $G$-map $f:M\rightarrow W$ in general
position with respect to $\Sigma $. The obstruction element $\mathfrak{o}%
_{G}^{n+1}(f|_{M^{n}})$ is the equivariant Poincar\'{e} dual of the point
set $f^{-1}\left( f(M)\cap \Sigma \right) $ understood as a chain in the
group
\begin{equation*}
\mathfrak{H}_{0}^{G}(M,H_{n}(W\backslash \Sigma ,%
\mathbb{Z}
)\otimes \mathcal{Z})
\end{equation*}%
with the appropriate coefficients in the group of coefficients.
\end{proposition}

\subparagraph{\textbf{4.}}

\noindent To compute the obstruction cocycle and the obstruction element, we
have to choose at least one equivariant cell structure on $M$ compatible
with the given action. We choose two equivariant cell structures connected
by a cellular map. To simplify the exposition, from now on we assume that $%
\mathcal{Z}$ is a trivial $G$-module.

\noindent Usually, the first equivariant cell structure induced on $M$ is a
simplicial one. It is used to define a piecewise affine $G$-map in general
position $f:M\rightarrow W$. The advantage of the simplicial structure is
that the map is completely determined by the images of the vertex orbits and
the requirement that the map is piecewise affine.

\noindent The second equivariant cell structure should satisfy the
requirement that the top dimensional group of chains is generated
equivariantly by a single cell $e$, and therefore
\begin{equation*}
C_{n+1}^{G}(M,%
\mathbb{Z}
)=%
\mathbb{Z}
\lbrack G]e.
\end{equation*}%
Such a structure will be called an \textbf{economic} $G$-structure of $M$.
The economic $G$-cell structure does not have to exist.

\subparagraph{5.}

\noindent The obstruction cocycle $\mathfrak{o}_{G}^{n+1}(f)$ is computed
using the simplicial cell structure and the geometric definition of the
obstruction cocycle.

\begin{proposition}
\label{Th:ObCocycle}Let $f:M\rightarrow W$ be a map in general position in
respect to $\Sigma $. Then for a $(n+1)$-simplex $\sigma $ of $M$ the
following formula holds
\begin{equation}
\mathfrak{o}_{G}^{n+1}(f)(\sigma )=\sum_{x\in f^{-1}(f(\sigma )\cap \Sigma )}%
\mathrm{I}(e,S_{f(x)})\left\Vert f(x),f(\sigma )\right\Vert \in
H_{n}(W\backslash \Sigma ,%
\mathbb{Z}
).  \label{ObCocycle}
\end{equation}%
Here $\mathrm{I}(e,S_{f(x)})$ denotes the intersection number of the image $%
f(e)$ and the appropriate oriented element $S_{f(x)}$ of the arrangement $%
\mathcal{S}$.
\end{proposition}

\noindent A cellular map between two cell structures allows evaluation of
the obstruction cocycle in both structures. In an economic $G$-structure,
every maximal cochain is determined by its value on the equivariant
generator $e$. Therefore the obstruction cocycle can be treated as an
element of the coefficient group $\mathfrak{o}_{G}^{n+1}(f)(e)\in
H_{n}(W\backslash \Sigma ,%
\mathbb{Z}
)$ expressed as a linear combination of the point classes. The obstruction
element is the class $[\mathfrak{o}_{G}^{n+1}(f)(e)]$, along the quotient
homomorphism $H_{n}(W\backslash \Sigma ,%
\mathbb{Z}
)\rightarrow H_{n}(W\backslash \Sigma ,%
\mathbb{Z}
)_{G}$.

\begin{theorem}
There exists a $G$-map $M\rightarrow W\backslash \Sigma $\ if and only if $%
[o_{G}^{n+1}(f)(e)]\in H_{n}(W\backslash \Sigma ,%
\mathbb{Z}
)_{G}$\ is zero.
\end{theorem}

\subsection{Example: Lov\'{a}sz conjecture.}

\noindent The first proof of the Lov\'{a}sz conjecture was given by Babson
and Kozlov \cite{Babson-Kozlov}. A revealingly simple proof was given by
Schultz, \cite{Carsten}, \cite{carsten-2}. Kozlov gave another proof in \cite%
{Kozlov}.

\begin{theorem}
Let $\Gamma $ be a graph, $r\geq 1$ and $n\geq 2$. If $\mathrm{Hom}%
(C_{2r+1},\Gamma )$ is $(n-1)$-connected, then $\Gamma $ is not $(n+2)$%
-colorable.
\end{theorem}

\noindent Here $C_{2r+1}$ is a circular graph with $2r+1$ edges. Let us
assume that $\Gamma $ is $(n+2)$-colorable. This means that there exists a
graph homomorphism $\Gamma \rightarrow K_{n+2}$, and consequently a map $%
\mathrm{Hom}(H,\Gamma )\rightarrow \mathrm{Hom}(H,K_{n+2})$ for every graph $%
H$. When we put $C_{2r+1}$with a $%
\mathbb{Z}
_{2}$-action instead of $H$, we end up with the $%
\mathbb{Z}
_{2}$-equivariant map%
\begin{equation}
\mathrm{Hom}(C_{2r+1},\Gamma )\rightarrow \mathrm{Hom}(C_{2r+1},K_{n+2})%
\text{.}  \label{map1}
\end{equation}%
The assumption that $\mathrm{Hom}(C_{2r+1},\Gamma )$ is $(n-1)$-connected
implies the existence of a $%
\mathbb{Z}
_{2}$ map%
\begin{equation}
S^{n}\rightarrow \mathrm{Hom}(C_{2r+1},\Gamma )\text{.}  \label{map2}
\end{equation}%
where $S^{n}$ is equipped with the antipodal action. Thus to prove the Lovas%
\'{z} conjecture it is enough to prove that there is no $%
\mathbb{Z}
_{2}$-map
\begin{equation}
S^{n}\rightarrow \mathrm{Hom}(C_{2r+1},K_{n+2}).  \label{map3}
\end{equation}%
In \cite{Carsten}, C. Schultz proved the nonexistence of the map (\ref{map3}%
) by comparing the complex $\mathrm{Hom}(C_{2r+1},K_{n+2})$ with the
complement of the torus arrangement and then performing some characteristic
class computations. Let us reproduce this beautiful construction and
substitute characteristic class computations with obstruction theory. Let $%
X_{r,n}=\left( S^{n}\right) ^{r}$ be a torus, and $A_{r,n}$ the union of the
arrangement of the following $r$ subtoruses%
\begin{equation*}
A_{r,n}^{i}=\{x\in X_{r,n}~|~x_{i}=-x_{i+1}\}\text{, }0\leq i\leq r-2\text{
and }A_{r,n}^{r-1}=\{x\in X_{r,n}~|~x_{r-1}=x_{0}\}\text{.}
\end{equation*}%
If we define a $%
\mathbb{Z}
_{2}=\langle \omega \rangle $ action on $X_{r,n}$ by%
\begin{equation*}
\omega \cdot (x_{0},x_{1},x_{2},..,x_{r-1})=(-x_{0},x_{r-1},x_{r-2},..,x_{1})
\end{equation*}%
it is apparent that $A_{r,n}$ is a $%
\mathbb{Z}
_{2}$-invariant subspace of $X_{r,n}$. Proposition 2.9 of \cite{Carsten}
says that there exists a $%
\mathbb{Z}
_{2}$-map%
\begin{equation}
\mathrm{Hom}(C_{2r+1},K_{n+2})\rightarrow X_{r,n}\backslash A_{r,n}\text{.}
\label{map4}
\end{equation}%
Thus, the assumption that $\Gamma $ is $(n+2)$-colorable implies the
existence of a $%
\mathbb{Z}
_{2}$-equivariant map%
\begin{equation*}
S^{n}\rightarrow X_{r,n}\backslash A_{r,n}.
\end{equation*}%
The following theorem provides a contradiction which implies Lov\'{a}sz
conjecture.

\begin{theorem}
There is no $%
\mathbb{Z}
_{2}$-equivariant map $S^{n}\rightarrow X_{r,n}\backslash A_{r,n}$.
\end{theorem}

\begin{proof}
The codimension of $A_{r,n}$ inside $X_{r,n}=\left( S^{n}\right) ^{r}$ is $n$%
. Since $A_{r,n}$ is compact, locally contractible space and $H^{i}(A_{r,n},%
\mathbb{Z}
)=0$ for $i>r(n-1)$ then Poincar\'{e} - Lefschetz Duality, (\cite[Corollary
8.4. p.352]{Bredon}, \cite[Proposition 3.46. p.254]{hatcher}), combined with
long exact sequence in cohomology, implies that the complement $%
X_{r,n}\backslash A_{r,n}$ is $n-2$ connected. Thus, the existence of a $%
\mathbb{Z}
_{2}$-map $S^{n}\rightarrow X_{r,n}\backslash A_{r,n}$ is determined by the
primary obstruction. Let $\xi =(0,..,0,1)\in S^{n}$. A map $%
f:S^{n}\rightarrow X_{r,n}$ defined by%
\begin{equation*}
f(x)=(x,\xi ,..,\xi )
\end{equation*}%
is a $%
\mathbb{Z}
_{2}$-map in general position and does not hit lower strata of the
arrangement. Indeed, the intersection contains two points%
\begin{equation*}
f(S^{n})\cap A_{r,n}=\{(\xi ,\xi ,..,\xi ),(-\xi ,\xi ,..,\xi )\}.
\end{equation*}%
Moreover,
\begin{equation*}
\begin{array}{lll}
(\xi ,\xi ,..,\xi )\in A_{r,n}^{r-1}\backslash \bigcup\limits_{i\neq
r-1}A_{r,n}^{i} & , & (-\xi ,\xi ,..,\xi )\in A_{r,n}^{0}\backslash
\bigcup\limits_{i\neq 0}A_{r,n}^{i}%
\end{array}%
\end{equation*}%
and $\omega \cdot (\xi ,..,\xi )=(-\xi ,\xi ,..,\xi )$. Proposition \ref%
{prop:SingularSet} implies that the obstruction element living in
\begin{equation*}
\mathfrak{H}_{%
\mathbb{Z}
_{2}}^{n}(S^{n},H_{n-1}(X_{r,n}\backslash A_{r,n},%
\mathbb{Z}
))
\end{equation*}%
is the equivariant Poincar\'{e} dual of the (orbit of the) point in
\begin{equation*}
\mathfrak{H}_{0}^{%
\mathbb{Z}
_{2}}(S^{n},H_{n-1}(X_{r,n}\backslash A_{r,n},%
\mathbb{Z}
)\otimes \mathcal{Z}).
\end{equation*}%
The isomorphisms \cite[Exercise 1, p.44]{Brown}, \cite[(1.5), p.57]{Brown},
and the fact that
\begin{equation*}
H_{n-1}(X_{r,n}\backslash A_{r,n},%
\mathbb{Z}
)\neq 0
\end{equation*}%
imply%
\begin{equation*}
\begin{array}{ll}
\mathfrak{H}_{0}^{%
\mathbb{Z}
_{2}}(S^{n},H_{n-1}(X_{r,n}\backslash A_{r,n},%
\mathbb{Z}
)\otimes \mathcal{Z}) & \cong H_{0}(%
\mathbb{Z}
_{2},H_{n-1}(X_{r,n}\backslash A_{r,n},%
\mathbb{Z}
)\otimes \mathcal{Z}) \\
& \cong \left( H_{n-1}(X_{r,n}\backslash A_{r,n},%
\mathbb{Z}
)\otimes \mathcal{Z}\right) _{%
\mathbb{Z}
_{2}} \\
& \neq 0.%
\end{array}%
\end{equation*}%
Since the orbit of a point is a generator of the $0$-homology, the
obstruction element is not zero and the map does not exist.
\end{proof}

\begin{remark}
In this example we were fortunate, because the intersection of the image of
the general position map and the forbidden set was just an orbit of a point.
That allowed a direct application of Proposition \ref{prop:SingularSet}.
\end{remark}

\section{\label{Result1}Partition of Sphere measures by Hyperplanes}

\noindent This chapter is contains an extension of Makeev's result \cite{Mak}%
.

\bigskip

\noindent A \textbf{proper} Borel measure on the sphere $S^{2}$ is a Borel
measure $\mu $ such that

(A) $\mu ([a,b])=0$ for any circular arc $[a,b]\subset S^{2}$, and

(B) $\mu (U)>0$ for each nonempty open set $U\subset S^{2}$.

\noindent Let $\Pi _{1}$, $\Pi _{2}$ and $\Pi _{3}$ be three planes through
the origin in $%
\mathbb{R}
^{3}$. Planes are in a \textbf{fan position} if they intersect along a
common line. The planes in a fan position cut the sphere $S^{2}$ in six
parts $\sigma _{1},..,\sigma _{6}$ naturally oriented up to a cyclic
permutation.

\begin{problem}
\label{Problem1}Find all the six-tuples $\mathbf{\alpha }=(\alpha
_{1},..,\alpha _{6})\in \mathbb{N}^{6}$ such that $\tfrac{\alpha _{1}}{%
\alpha _{1}+...+\alpha _{6}}+...+\tfrac{\alpha _{6}}{\alpha _{1}+...+\alpha
_{6}}=1$, and that for any proper Borel probability measure $\mu $ on the
sphere $S^{2}$ there exist three planes in a fan position, with angular
sectors having the prescribed amount of the measure, i.e., for all $i{\small %
\in \{1,..,6\}}$,
\begin{equation*}
\mu (\sigma _{i})=\tfrac{\alpha _{i}}{\alpha _{1}+...+\alpha _{6}}.
\end{equation*}%
The six-tuples which satisfy these conditions are \textbf{solutions} of the
problem.
\end{problem}

\noindent The existence of the equipartition solution was proved by V. V.
Makeev \cite{Mak}. Modifying the configuration space and the test map, we
prove in section \ref{Sec:ThMain1} that beside $(1,1,1,1,1,1)$ there exist
at least two more solutions.

\begin{theorem}
\label{th:Main1}Let $\mu $ be a \textit{proper} Borel probability measure on
the sphere $S^{2}$. Then there are three planes intersecting along a line
such that the ratio of the measure $\mu $ in the angular sectors cut by the
planes is%
\begin{equation*}
(A)~(1,1,2,1,1,2)\qquad \qquad (B)~(1,2,2,1,2,2).
\end{equation*}
\end{theorem}

\subsection{The configuration space / test map scheme}

\noindent We use the configuration space / test map scheme to reduce the
partition problem to an equivariant one. The basic idea, which we modify,
comes from the papers of Imre B\'{a}r\'{a}ny and Ji\v{r}i Matou\v{s}ek \cite%
{BaMa2001}, \cite{BaMa2002}.

\bigskip

\noindent The $k$\textbf{-fan} $(l;\Pi _{1},\Pi _{2},\ldots ,\Pi _{k})$ in $%
\mathbb{R}
^{3}$ (or on $S^{2}$ ) is formed from an oriented line through the origin $l$
and $k$ closed half planes $\Pi _{1},\Pi _{2},\ldots ,\Pi _{k}$ which
intersect along the common boundary $l=\partial \Pi _{1}=\ldots =\partial
\Pi _{k}$. The intersection of a $k$-fan with the sphere $S^{2}$ is also
equally called a $k$-fan. Thus, the collection $(x;l_{1},\ldots ,l_{k})$ of
a point $x\in S^{2}$ and $k$ great semicircles $l_{1},\ldots ,l_{k}$
emanating from $x$ is also a $k$-fan. Sometimes instead of great semicircles
we use:

(A) open angular sectors $\sigma _{i}$ between $l_{i}$ and $%
l_{i+1},\,i=1,\ldots ,k$; or

(B) tangent vectors $t_{i}\in T_{x}S^{2}$ which are determined by the great
semicircle curves $l_{i}$, $\,i=1,\ldots ,k$.

\noindent Here $T_{x}S^{2}$ denotes the tangent space at a point $x\in S^{2}$%
. So there are four equally useful notations for a single $k$-fan, and we
prefer the tangent vector notation $(x;t_{1},\ldots ,t_{k})$. The space of
all $k$-fans in $%
\mathbb{R}
^{3}$ will be denoted by $F_{k}$.

\bigskip

\noindent \textbf{The configuration space.} For a proper Borel probability
measure $\mu $ on $S^{2}$ and $n>1$, the configuration space is defined by
\begin{equation*}
X_{\mu ,n}=\{(x;t_{1},\ldots ,t_{n})\in F_{n}\mid (\forall i=1,\ldots
,n)\,\mu (\sigma _{i})=\tfrac{{\small 1}}{{\small n}}\}.
\end{equation*}%
Since every $n$-fan $(x;t_{1},\ldots ,t_{n})$ of the configuration space $%
X_{\mu ,n}$ is completely determined by the pair $(x,t_{1})\in S^{2}\times
T_{x}S^{2}$ and the measure $\mu $, there is a homeomorphism $X_{\mu
,n}\cong V_{2}(%
\mathbb{R}
^{3})$.

\noindent \textbf{The test map.} Let us fix the "symmetric" six-tuple $%
\mathbf{\alpha }=(\alpha _{1},\alpha _{2},\alpha _{3},\alpha _{1},\alpha
_{2},\alpha _{3})\in \mathbb{N}^{6}$ such that $\alpha _{1}+\alpha
_{2}+\alpha _{3}=\frac{n}{2}$. The test map for our problem is defined by
\begin{equation*}
\begin{array}{ll}
\Phi :X_{\mu ,n}\rightarrow W_{n}=\{x\in
\mathbb{R}
^{n}\mid \Sigma _{i}x_{i}=0\}, & \Phi (x;t_{1},\ldots ,t_{n})=(\theta _{1}-%
\tfrac{2\pi }{{\small n}},\ldots ,\theta _{n}-\tfrac{2\pi }{{\small n}}),%
\end{array}%
\end{equation*}%
where $\theta _{i}$ is the angle between tangent vectors $t_{i}$ and $%
t_{i+1} $ in the tangent plane $T_{x}S^{2}$. Here we assume that $%
t_{n+1}=t_{1}$.

\noindent \textbf{The action.} The dihedral group $\mathbb{D}_{2n}=\langle
j,\varepsilon \,|\,\varepsilon ^{n}=j^{2}=1,\,\varepsilon j=j\varepsilon
^{n-1}\,\rangle $ acts both on the configuration space $X_{\mu ,n}$ and on
the hyperplane $W_{n}$ in the following way
\begin{equation*}
\left\{
\begin{array}{l}
\varepsilon (x;t_{1},\ldots ,t_{n})=(x;t_{n},t_{1},\ldots ,t_{n-1}) \\
j(x;t_{1},\ldots ,t_{n})=(-x;t_{1},t_{n},\ldots ,t_{2})%
\end{array}%
,\right. \left\{
\begin{array}{c}
\varepsilon (x_{1},\ldots ,x_{n})=(x_{n},x_{1}\ldots ,x_{n-1}) \\
j(x_{1},\ldots ,x_{n})=(x_{n},\ldots ,x_{2},x_{1})%
\end{array}%
,\right.
\end{equation*}%
for $(x;t_{1},\ldots ,t_{n})\in X_{\mu ,n}$ and $(x_{1},\ldots ,x_{n})\in
W_{n}$. The action of $\mathbb{D}_{2n}$ on $X_{\mu ,n}$ is \textbf{free} and
the test map $\Phi $ is equivariant.

\noindent \textbf{The test space. }The test space in this symmetric problem
is the union $\bigcup \mathcal{A}\subset W_{n}$ of the smallest $\mathbb{D}%
_{2n}$-invariant arrangement $\mathcal{A}$, which contains a linear subspace
$L\subset W_{n}$ defined by the equalities
\begin{equation*}
\begin{array}{lll}
x_{1}+\ldots +x_{\frac{n}{2}}= & x_{\alpha _{1}+1}+\ldots +x_{\alpha _{1}+%
\frac{n}{2}}= & x_{\alpha _{1}+\alpha _{2}+1}+\ldots +x_{\alpha _{1}+\alpha
_{2}+\frac{n}{2}}=0\text{.}%
\end{array}%
\end{equation*}

\noindent We have proved the basic proposition of the configuration space /
test map scheme.

\begin{proposition}
\label{prop:VezaProblem-Ekvivarijantan}Let $\mathbf{\alpha }=(\alpha
_{1},\alpha _{2},\alpha _{3},\alpha _{1},\alpha _{2},\alpha _{3})\in \mathbb{%
N}^{6}$ be a symmetric $6$-tuple such that $\alpha _{1}+\alpha _{2}+\alpha
_{3}=\frac{n}{2}$. If there is no $\mathbb{D}_{2n}$-equivariant map%
\begin{equation*}
V_{2}(%
\mathbb{R}
^{3})\rightarrow W_{n}\setminus \bigcup \mathcal{A,}
\end{equation*}%
then for every proper Borel probability measure on the sphere $S^{2}$ there
exist three planes in a fan position with angular sectors such that%
\begin{equation*}
(\forall i{\small \in \{1,..,6\}})~\mu (\sigma _{i})=\tfrac{\alpha _{i}}{n}.
\end{equation*}
\end{proposition}

\noindent The extension of scalars from homological algebra (as shown in
\cite{Bl} and \cite{Bl-Vr-Ziv}) allows us to prove the following
equivalence. Let $\mathbb{Q}_{4n}$ denote the generalized quaternion group $%
\langle \epsilon ,j\rangle \subset S^{3}$, (section \ref{Sec:Q4n}).

\begin{proposition}
\label{prop:PrelazakNaNovuGrupu}Following maps jointly exist or do not exist:%
\begin{equation*}
\text{a }\mathbb{D}_{2n}\text{-map }V_{2}(%
\mathbb{R}
^{3})\rightarrow W_{n}\setminus \bigcup \mathcal{A}\text{ \ and a }\mathbb{Q}%
_{4n}\text{-map }S^{3}\rightarrow W_{n}\setminus \bigcup \mathcal{A}.
\end{equation*}%
The group $\mathbb{Q}_{4n}$ acts on $S^{3}$ as a subgroup and on $W_{n}$ via
the quotient homomorphism $\mathbb{Q}_{4n}\rightarrow \mathbb{Q}_{4n}/H\cong
\mathbb{D}_{2n}$ (see Appendix).
\end{proposition}

\subsection{\label{Sec:ThMain1}Proof of Theorem \protect\ref{th:Main1}}

\noindent According to Propositions \ref{prop:VezaProblem-Ekvivarijantan}
and \ref{prop:PrelazakNaNovuGrupu} it is enough to prove that there is no $%
\mathbb{Q}_{4n}$-map $S^{3}\rightarrow W_{n}\setminus \bigcup \mathcal{A}$.
Here $\mathcal{A}$ is the minimal $\mathbb{Q}_{4n}$(=$\mathbb{D}_{2n}$)
arrangement containing the subspace $L$ defined

\noindent (A) for $n=8$ and $\mathbf{\alpha }=(\alpha _{1},\alpha
_{2},\alpha _{3},\alpha _{1},\alpha _{2},\alpha _{3})=(1,1,2,1,1,2)$, by the
equations:%
\begin{equation*}
\begin{array}{lll}
x_{1}+x_{2}+x_{3}+x_{4}~= & x_{2}+x_{3}+x_{4}+x_{5}~= &
x_{3}+x_{4}+x_{5}+x_{6}=0\text{;}%
\end{array}%
\end{equation*}

\noindent (B) for $n=10$ and $\mathbf{\alpha }=(\alpha _{1},\alpha
_{2},\alpha _{3},\alpha _{1},\alpha _{2},\alpha _{3})=(1,2,2,1,2,2)$, by the
equations:%
\begin{equation*}
\begin{array}{lll}
x_{1}+...+x_{5}~= & x_{2}+...+x_{6}~= & x_{4}+...+x_{8}~=0\text{.}%
\end{array}%
\end{equation*}%
The codimension of $\mathcal{A}$ inside $W_{n}$ in both cases is $3$, so the
complement $W_{n}\setminus \bigcup \mathcal{A}$ is $1$-connected (by Poincar%
\'{e} - Lefschetz Duality \cite[Corollary 8.4, p.352]{Bredon}, \cite[%
Proposition 3.46, p.254]{hatcher} or the Goresky-MacPherson formula).
Therefore the primary obstruction is responsible for the existence of $%
\mathbb{Q}_{4n}$-maps $S^{3}\rightarrow W_{n}\setminus \bigcup \mathcal{A}$,
and the obstruction exact sequence has the form%
\begin{equation*}
\lbrack S^{3},W_{n}\setminus \bigcup \mathcal{A}]_{\mathbb{Q}%
_{4n}}\longrightarrow \{\ast \}\overset{\mathfrak{o}_{\mathbb{Q}_{4n}}^{3}}{%
\longrightarrow }\mathfrak{H}_{\mathbb{Q}_{4n}}^{3}(S^{3},H_{2}\left(
W_{n}\setminus \bigcup \mathcal{A},%
\mathbb{Z}
\right) ).
\end{equation*}%
Since $\mathbb{Q}_{4n}\subset S^{3}$ acts on $S^{3}$ as a subgroup and $%
S^{3} $ is connected, the $\mathbb{Q}_{4n}$-module $\mathcal{Z}$ is trivial.
The Equivariant Poincar\'{e} duality isomorphism (Theorem \ref{Th:EqPoincare}%
) implies%
\begin{equation*}
\begin{array}{ll}
\mathfrak{H}_{\mathbb{Q}_{4n}}^{3}(S^{3},H_{2}\left( W_{n}\setminus \bigcup
\mathcal{A},%
\mathbb{Z}
\right) ) & \cong \mathfrak{H}_{0}^{\mathbb{Q}_{4n}}(S^{3},H_{2}\left(
W_{n}\setminus \bigcup \mathcal{A},%
\mathbb{Z}
\right) \otimes \mathcal{Z}) \\
& \cong H_{2}\left( W_{n}\setminus \bigcup \mathcal{A},%
\mathbb{Z}
\right) _{\mathbb{Q}_{4n}}\text{.}%
\end{array}%
\end{equation*}%
Thus, in both cases we are going to define a map in general position,
compute the obstruction cocycle and identify its obstruction element inside
a group of coinvariants. All computations are done by the Mathematica 5.0
package.

\subsubsection{Case $n=8$ and $\mathbf{\protect\alpha }=(1,1,2,1,1,2)$}

Let us define a map $f:S^{3}\rightarrow W_{8}$ on the vertex $t$ ($=a_{1}$
in Appendix) by
\begin{equation*}
f(t)=(-3,3,-1,1,1,-2,2,-1)
\end{equation*}%
and extend it equivariantly. For example $f(jt)=(-1,2,-2,1,1,-1,3,-3)$.

\begin{figure}[tbh]
\centering\includegraphics[scale=0.80]{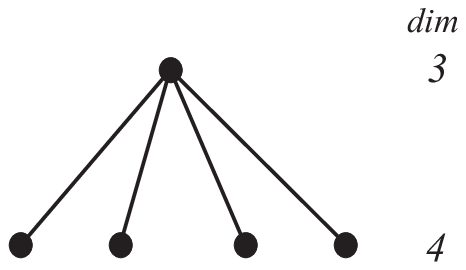} \caption{{}}
\label{fig:Fig22}
\end{figure}

\noindent For the subspace $L$ defined by

\begin{equation*}
x_{1}+x_{2}+x_{3}+x_{4}=x_{2}+x_{3}+x_{4}+x_{5}=x_{3}+x_{4}+x_{5}+x_{6}=%
\sum_{i=1}^{8}x_{i}=0\text{,}
\end{equation*}%
the arrangement $\mathcal{A}$ is the minimal $\mathbb{Q}_{32}$-arrangement
containing $L$. It has four maximal elements $L$, $\epsilon L$, $\epsilon
^{2}L$ and $\epsilon ^{3}L$. This can easily be seen from the set equalities
$\epsilon ^{4}L=L$ and $\epsilon ^{2}L=jL$. The intersection $L\cap \epsilon
L\cap \epsilon ^{2}L\cap \epsilon ^{3}L$ is a codimension $1$ subspace of $L$%
, $\epsilon L$, $\epsilon ^{2}L$, $\epsilon ^{3}L$. Thus, the Hasse diagram
of the intersection poset of the arrangement $\mathcal{A}$ with the reversed
order is as in Figure \ref{fig:Fig22}.

\noindent Now we intersect the image under $f$ of the maximal cell $e=\left(
[t,\epsilon t]\cup ..\cup \lbrack \epsilon ^{7}t,\epsilon ^{8}t]\right) \ast
\lbrack jt,\epsilon jt]$ with the test space $\bigcup \mathcal{A}=L\cup
\epsilon L\cup \epsilon ^{2}L\cup \epsilon ^{3}L$.\ The results of the $%
8\times 4=32$ intersections $f\left( [\epsilon ^{i}t,\epsilon ^{i+1}t]\ast
\lbrack jt,\epsilon jt]\right) \cap \epsilon ^{r}L$ are summarized in the
following table:

\begin{equation*}
\begin{tabular}{|l|l|l|l|}
\hline
& \multicolumn{3}{|l|}{\qquad \qquad\ $f\left( [\epsilon ^{i}t,\epsilon
^{i+1}t]\ast \lbrack jt,\epsilon jt]\right) \cap \epsilon ^{r}L$} \\ \hline
$\mathbf{i/r}$ & {\small preimage of the intersection point} & {\small %
intersection point in }$W_{8}$ &  \\ \hline
$1/2$ & $\frac{3}{7}jt+\frac{1}{14}\epsilon jt+\frac{1}{14}\epsilon t+\frac{3%
}{7}\epsilon ^{2}t$ & $\left( -\tfrac{1}{2},\frac{15}{14},-\frac{2}{7},-%
\frac{2}{7},\frac{15}{14},-\tfrac{1}{2},-\frac{2}{7},-\frac{2}{7}\right) $ &
$p_{1}$ \\ \hline
$2/2$ & $\frac{20}{51}jt+\frac{23}{153}\epsilon jt+\frac{62}{153}\epsilon
^{2}t+\frac{8}{153}\epsilon ^{3}t$ & $\left( -\tfrac{4}{9},\tfrac{16}{17},-%
\frac{1}{3},-\frac{25}{153},1,-\tfrac{77}{153},-\frac{1}{3},-\tfrac{25}{153}%
\right) $ & $p_{2}$ \\ \hline
$2/1$ & $\frac{8}{153}jt+\frac{63}{153}\epsilon jt+\frac{23}{153}\epsilon
^{2}t+\frac{20}{51}\epsilon ^{3}t$ & $\left( 1,-\frac{25}{153},-\frac{1}{3},%
\frac{16}{17},-\frac{4}{9},-\frac{25}{153},-\frac{1}{3},-\frac{77}{153}%
\right) $ & $p_{3}$ \\ \hline
$4/1$ & $\frac{1}{20}jt+\frac{11}{40}\epsilon jt+\frac{1}{10}\epsilon ^{4}t+%
\frac{23}{40}\epsilon ^{5}t$ & $\left( -\frac{11}{20},\frac{1}{2},-\frac{1}{5%
},-\frac{3}{2},\frac{6}{5},\frac{1}{2},-\frac{1}{5},\frac{1}{4}\right) $ & $%
p_{4}$ \\ \hline
$4/0$ & $\frac{23}{40}jt+\frac{1}{10}\epsilon jt+\frac{11}{40}\epsilon ^{4}t+%
\frac{1}{20}\epsilon ^{5}t$ & $\left( -\frac{1}{5},\frac{1}{2},-\frac{11}{20}%
,\frac{1}{4},-\frac{1}{5},\frac{1}{2},\frac{6}{5},-\frac{3}{2}\right) $ & $%
p_{5}$ \\ \hline
$5/0$ & $\frac{1}{3}jt+\frac{1}{6}\epsilon jt+\frac{1}{6}\epsilon ^{5}t+%
\frac{1}{3}\epsilon ^{6}t$ & $\left( \frac{1}{3},\frac{1}{3},-\frac{5}{3},1,%
\frac{1}{3},\frac{1}{3},1,-\frac{5}{3}\right) $ & $p_{6}$ \\ \hline
$7/3$ & $\frac{1}{14}jt+\frac{3}{7}\epsilon jt+\frac{3}{7}\epsilon ^{7}t+%
\frac{1}{14}t$ & $\left( \frac{1}{7},-\frac{25}{14},\frac{3}{2},\frac{1}{7},%
\frac{1}{7},\frac{3}{2},-\frac{25}{14},\frac{1}{7}\right) $ & $p_{7}$ \\
\hline
\end{tabular}%
\end{equation*}

\noindent The formula for the obstruction cocycle (\ref{ObCocycle}) implies
that there are signs $s_{1},..,s_{7}\in \{1,-1\}$ such that%
\begin{equation}
\mathfrak{o}_{\mathbb{Q}_{32}}(f)(e)=\sum_{i=1}^{7}s_{i}\left\Vert
p_{i}\right\Vert .  \label{CoCycleProblem1}
\end{equation}%
The point classes in (\ref{CoCycleProblem1}) do not depend on the embedding
of the associated simplices because each intersection point $p_{i}$ is
contained in just one element of the arrangement. In this situation, in
contrast to the Lovas\'{z} conjecture, we can not just apply Proposition \ref%
{prop:SingularSet}.

\noindent To prove that the obstruction element does or does not vanish, the
cohomology class of the obstruction cocycle inside $H_{2}(W_{n}\setminus
\bigcup \mathcal{A};\mathbb{%
\mathbb{Z}
})_{\mathbb{Q}_{4n}}$ has to be computed. However, to prove that the
obstruction element does not vanish, we are not compelled to completely
identify the obstruction element.

\bigskip

\noindent\ With a help of the Poincar\'{e}-Alexander duality isomorphism and
the Universal coefficient isomorphism we have
\begin{equation}
\begin{array}{ll}
H_{2}(W_{n}\setminus \bigcup \mathcal{A},\mathbb{%
\mathbb{Z}
}) & \cong H^{n-4}(\bigcup \widehat{\mathcal{A}},\mathbb{%
\mathbb{Z}
}) \\
& \cong \mathrm{Hom}\left( H_{n-4}\left( \bigcup \widehat{\mathcal{A}},%
\mathbb{%
\mathbb{Z}
}\right) ;\mathbb{%
\mathbb{Z}
}\right) \oplus \mathrm{Ext}\left( H_{n-5}\left( \bigcup \widehat{\mathcal{A}%
}\right) ,\mathbb{%
\mathbb{Z}
}\right)%
\end{array}
\label{PoincareDuality}
\end{equation}%
where $\widehat{\mathcal{A}}$ denotes the one-point compactification of the
arrangement $\mathcal{A}$. The calculations of $H_{n-4}\left( \bigcup
\widehat{\mathcal{A}};\mathbb{%
\mathbb{Z}
}\right) $ and $\mathrm{Ext}\left( H_{n-5}\left( \bigcup \widehat{\mathcal{A}%
}\right) ,\mathbb{%
\mathbb{Z}
}\right) $ can be carried out by the Ziegler-\v{Z}ivaljevi\'{c} formula \cite%
{ZZ}. For example, there is a decomposition
\begin{equation}
\begin{array}{ll}
H_{n-4}(\bigcup \widehat{\mathcal{A}},\mathbb{%
\mathbb{Z}
}) & \cong \underset{V\in P}{\bigoplus }H_{n-4}\left( \Delta (P_{<V})\ast
S^{\dim V},\mathbb{%
\mathbb{Z}
}\right) \\
& \cong \underset{d=0}{\overset{n-4}{\bigoplus }}\left( \underset{V\in
P:\dim V=d}{\bigoplus }\tilde{H}_{n-5-d}(\Delta (P_{<V}),\mathbb{%
\mathbb{Z}
}\mathbf{)}\right)%
\end{array}
\label{Z-z-decomposition}
\end{equation}%
where $P$ is the intersection poset of the arrangement $\mathcal{A}$. By
convention, $\tilde{H}_{-1}(\emptyset )=\mathbb{%
\mathbb{Z}
}$. The decomposition (\ref{Z-z-decomposition}) is not a decomposition of $%
\mathbb{Q}_{4n}$-modules. This fact is illustrated in \cite[Theorem 4.7.(C).]%
{Bl-Vr-Ziv}. If considered with the appropriate field coefficients (\ref%
{Z-z-decomposition}) is a decomposition of $\mathbb{Q}_{4n}$-modules \cite[%
Proposition 2.3]{Sund}. In order to use (\ref{PoincareDuality}) and (\ref%
{Z-z-decomposition}) and to compute the coinvariants $H_{2}(W_{n}\setminus
\bigcup \mathcal{A};\mathbb{%
\mathbb{Z}
})_{\mathbb{Q}_{4n}}$, we have to keep in mind that Poincar\'{e}-Alexander
duality isomorphism is not an isomorphism of $\mathbb{Q}_{4n}$-modules.
Fortunately, it is a $\mathbb{Q}_{4n}$-map up to an orientation character.
On the other hand, the universal coefficient isomorphism, when the $\mathrm{%
Ext}$ part vanishes, is a $\mathbb{Q}_{4n}$-map, if the action on $\mathrm{%
Hom}$ is given by%
\begin{equation*}
(g\cdot f)(x)=f(g^{-1}\cdot x)
\end{equation*}%
for $g\in \mathbb{Q}_{4n}$, $f\in \mathrm{Hom}\left( H_{n-4}\left( \bigcup
\widehat{\mathcal{A}},\mathbb{%
\mathbb{Z}
}\right) ;\mathbb{%
\mathbb{Z}
}\right) $ and $x\in H_{n-4}\left( \bigcup \widehat{\mathcal{A}},\mathbb{%
\mathbb{Z}
}\right) $.

\noindent Let a modified $\mathbb{Q}_{4n}$ action $\ast $ on $\mathrm{Hom}%
\left( H_{n-4}\left( \bigcup \widehat{\mathcal{A}},\mathbb{%
\mathbb{Z}
}\right) ;\mathbb{%
\mathbb{Z}
}\right) $ be defined by%
\begin{equation}
(g\ast f)(x)=\det (g)f(g^{-1}\cdot x)  \label{eq:modifiedAction}
\end{equation}%
for $g\in \mathbb{Q}_{4n}$, $f\in \mathrm{Hom}\left( H_{n-4}\left( \bigcup
\widehat{\mathcal{A}},\mathbb{%
\mathbb{Z}
}\right) ;\mathbb{%
\mathbb{Z}
}\right) $ and $x\in H_{n-4}\left( \bigcup \widehat{\mathcal{A}},\mathbb{%
\mathbb{Z}
}\right) $. If the $\mathrm{Ext}$ part vanishes in (\ref{PoincareDuality})
and the action on the $\mathrm{Hom}$ part is assumed to be $\ast $, the
Poincar\'{e} duality isomorphism
\begin{equation*}
H_{2}(W_{n}\setminus \bigcup \mathcal{A},\mathbb{%
\mathbb{Z}
})\cong \mathrm{Hom}\left( H_{n-4}\left( \bigcup \widehat{\mathcal{A}},%
\mathbb{%
\mathbb{Z}
}\right) ;\mathbb{%
\mathbb{Z}
}\right)
\end{equation*}%
becomes an isomorphism of $\mathbb{Q}_{4n}$-modules. Consequently,
\begin{equation*}
H_{2}(W_{n}\setminus \bigcup \mathcal{A};\mathbb{%
\mathbb{Z}
})_{\mathbb{Q}_{4n}}\cong \mathrm{Hom}\left( H_{n-4}\left( \bigcup \widehat{%
\mathcal{A}},\mathbb{%
\mathbb{Z}
}\right) ;\mathbb{%
\mathbb{Z}
}\right) _{\mathbb{Q}_{4n}}.
\end{equation*}

\bigskip

\noindent Now we identify every point class from the sum$\mathfrak{\ }$(\ref%
{CoCycleProblem1}) using the isomorphism (\ref{PoincareDuality})
\begin{equation*}
\varphi :H_{2}\left( W_{n}\setminus \bigcup \mathcal{A},\mathbb{%
\mathbb{Z}
}\right) \rightarrow \mathrm{Hom}\left( H_{n-4}\left( \bigcup \widehat{%
\mathcal{A}},\mathbb{%
\mathbb{Z}
}\right) ;\mathbb{%
\mathbb{Z}
}\right) .
\end{equation*}%
The isomorphism is an evaluation of the linking number (when it is correctly
defined). For point classes the following formula holds
\begin{equation*}
\varphi \left( \left\Vert f(x)\right\Vert \right) (l)=\mathrm{link}%
(l,\left\Vert f(x)\right\Vert )l
\end{equation*}%
where $l\in H_{n-4}(\bigcup \widehat{\mathcal{A}};\mathbb{%
\mathbb{Z}
})$. Finally we must determine whether the class of
\begin{equation*}
\varphi \left( \mathfrak{o}_{\mathbb{Q}_{4n}}(f)(e)\right) \in \mathrm{Hom}%
\left( H_{n-4}\left( \bigcup \widehat{\mathcal{A}},\mathbb{%
\mathbb{Z}
}\right) ;\mathbb{%
\mathbb{Z}
}\right)
\end{equation*}%
in the group of coinvariants $\mathrm{Hom}\left( H_{n-4}\left( \bigcup
\widehat{\mathcal{A}},\mathbb{%
\mathbb{Z}
}\right) ;\mathbb{%
\mathbb{Z}
}\right) _{\mathbb{Q}_{4n}}$ is or is not zero. In each case this is a
different, and usually a very difficult problem.

\bigskip

\noindent In the present situation, the exact sequence of $\mathbb{Q}_{4n}$%
-modules%
\begin{equation*}
0\rightarrow \underset{V\in P:\dim V=n-4}{\bigoplus }\tilde{H}%
_{n-5-d}(\Delta (P_{<V});\mathbb{%
\mathbb{Z}
}\mathbf{)}\mathbb{\rightarrow }H_{n-4}\left( \bigcup \widehat{\mathcal{A}},%
\mathbb{%
\mathbb{Z}
}\right)
\end{equation*}%
induces an exact sequence of $\mathbb{Q}_{4n}$-modules
\begin{equation}
H_{2}\left( W_{n}\setminus \bigcup \mathcal{A},\mathbb{%
\mathbb{Z}
}\right) \overset{\chi }{\rightarrow }\mathrm{Hom}\left( \underset{V\in
P:\dim V=n-4}{\bigoplus }\tilde{H}_{n-5-d}(\Delta (P_{<V});\mathbb{%
\mathbb{Z}
}\mathbf{)},\mathbb{%
\mathbb{Z}
}\right) \rightarrow 0\text{.}  \label{HomPartProjection}
\end{equation}%
The geometric interpretation of the map $\chi $ is the computation of the
linking numbers. Let us assume that all of the $(n-5-d)$ homology groups in
the above sum are free. The left exactness of the coinvariant functor
implies that the sequence%
\begin{equation}
H_{2}\left( W_{n}\setminus \bigcup \mathcal{A},\mathbb{%
\mathbb{Z}
}\right) _{\mathbb{Q}_{4n}}\overset{\chi ^{\ast }}{\rightarrow }\mathrm{Hom}%
\left( \underset{V\in P:\dim V=n-4}{\bigoplus }\tilde{H}_{n-5-d}(\Delta
(P_{<V});\mathbb{%
\mathbb{Z}
}\mathbf{)},\mathbb{%
\mathbb{Z}
}\right) _{\mathbb{Q}_{4n}}~\rightarrow ~0  \label{projection}
\end{equation}%
is exact. The geometric interpretation of the map $\chi ^{\ast }$ is the
summation of linking. The map $\chi ^{\ast }$ is a good test map for
detecting whether an element $\mathfrak{o}\in H_{2}\left( W_{n}\setminus
\bigcup \mathcal{A},\mathbb{%
\mathbb{Z}
}\right) _{\mathbb{Q}_{4n}}$ is not zero.

\noindent We use the map $\chi ^{\ast }$ (\ref{projection}) and prove that $%
\chi ^{\ast }([\mathfrak{o}_{\mathbb{Q}_{32}}(f)])\neq 0$. Let us identify
groups in question (\ref{HomPartProjection}) and (\ref{projection}).

\begin{lemma}
(A)$\underset{V\in P:\dim V=4}{\bigoplus }\tilde{H}_{-1}(\Delta (P_{<V});%
\mathbb{%
\mathbb{Z}
}\mathbf{)~}\mathbb{\cong ~%
\mathbb{Z}
}^{4}$, as abelian groups.

(B) $\mathrm{Hom}\left( \underset{V\in P:\dim V=4}{\bigoplus }\tilde{H}%
_{-1}(\Delta (P_{<V});\mathbb{%
\mathbb{Z}
}\mathbf{)},\mathbb{%
\mathbb{Z}
}\right) _{\mathbb{Q}_{32}}~\mathbb{\cong ~%
\mathbb{Z}
}_{2}$.
\end{lemma}

\begin{proof}
Statement (A) follows from the Hasse diagram of the arrangement and the
Ziegler-\v{Z}ivaljevi\'{c} formula. Statement (B) is a consequence of the
set equality $\epsilon ^{4}L=L$ and the following orientation computation:
the element $\epsilon ^{4}$ acts on $W_{8}$ by preserving its orientation.
On the orthogonal complement $L^{\perp }$ of $L$ the operator $\epsilon ^{4}$%
, for the basis $%
\{e_{1}+e_{2}+e_{3}+e_{4},~e_{2}+e_{3}+e_{4}+e_{5},~e_{3}+e_{4}+e_{5}+e_{6},~e_{1}+...+e_{8}\}
$ of $L^{\perp }$, has the matrix%
\begin{equation*}
\Xi =\left(
\begin{array}{cccc}
{\small -1} & {\small 0} & {\small 0} & {\small 1} \\
{\small 0} & {\small -1} & {\small 0} & {\small 1} \\
{\small 0} & {\small 0} & {\small -1} & {\small 1} \\
{\small 0} & {\small 0} & {\small 0} & {\small 1}%
\end{array}%
\right) \text{.}
\end{equation*}%
Since $\det \Xi =-1$, the element $\epsilon ^{4}$ changes the orientation of
$L^{\perp }$ and consequently it changes the orientation on $L$. If $l\in
\underset{V\in P:\dim V=4}{\bigoplus }\tilde{H}_{-1}(\Delta (P_{<V});\mathbb{%
\mathbb{Z}
}\mathbf{)}$ is the generator associated with the subspace $L$, then the sum
$\underset{V\in P:\dim V=4}{\bigoplus }\tilde{H}_{-1}(\Delta (P_{<V});%
\mathbb{%
\mathbb{Z}
}\mathbf{)}$ is generated as a free abelian group by elements $l$, $\epsilon
l$, $\epsilon ^{2}l$, $\epsilon ^{3}l$. The set equality $\epsilon ^{4}L=L$
implies the homology equality $\epsilon ^{4}l=-l$.

\noindent Let $\xi \in $ $\mathrm{Hom}\left( \underset{V\in P:\dim V=4}{%
\bigoplus }\tilde{H}_{-1}(\Delta (P_{<V});\mathbb{%
\mathbb{Z}
}\mathbf{)},\mathbb{%
\mathbb{Z}
}\right) $ be given by $\xi (l)=1$ and $\xi (\epsilon ^{i}l)=0$ for $i\in
\{1,2,3\}$. Then by (\ref{eq:modifiedAction})
\begin{equation*}
\epsilon ^{4}\ast \xi (l)=\det (\epsilon ^{4})\xi (\epsilon ^{-4}l)=\xi
(\epsilon ^{4}l)=\xi (-l)=-1
\end{equation*}%
and
\begin{equation*}
\epsilon ^{4}\ast \xi (\epsilon ^{i}l)=\det (\epsilon ^{4})\xi (\epsilon
^{-4+i}l)=0
\end{equation*}%
for $i\in \{1,2,3\}$. Thus, there is a relation
\begin{equation*}
\xi +\epsilon ^{4}\ast \xi =0
\end{equation*}%
which in coinvariants becomes $2[\xi ]=0$. The statement (B) follows.
\end{proof}

\noindent The proof of case (A) of Theorem \ref{th:Main1} is a consequence
of the following lemma.

\begin{lemma}
(A) The element $\chi (\mathfrak{o}_{\mathbb{Q}_{32}}(f))$ is given by%
\begin{equation*}
\begin{array}{ll}
\chi (\mathfrak{o}_{\mathbb{Q}_{32}}(f))(l)=s_{5}+s_{6}, & \chi (\mathfrak{o}%
_{\mathbb{Q}_{32}}(f))(\epsilon ^{1}l)=s_{3}+s_{4}, \\
\chi (\mathfrak{o}_{\mathbb{Q}_{32}}(f))(\epsilon ^{2}l)=s_{1}+s_{2}, & \chi
(\mathfrak{o}_{\mathbb{Q}_{32}}(f))(\epsilon ^{3}l)=1.%
\end{array}%
\end{equation*}

(B) $\chi ^{\ast }([\mathfrak{o}_{\mathbb{Q}_{32}}(f)])\in \mathrm{Hom}%
\left( \underset{V\in P:\dim V=4}{\bigoplus }\tilde{H}_{-1}(\Delta (P_{<V});%
\mathbb{%
\mathbb{Z}
}\mathbf{)},\mathbb{%
\mathbb{Z}
}\right) _{\mathbb{Q}_{32}}$ is the generator of the group $\mathbb{%
\mathbb{Z}
}_{2}$.
\end{lemma}

\begin{proof}
Statement (A) follows from formula (\ref{CoCycleProblem1}) and the
associated table of intersections. The second statement is a consequence of
the geometric interpretation of the map $\chi ^{\ast }$ - summation of the
intersection numbers.
\end{proof}

\noindent Thus the obstruction element $[\mathfrak{o}_{\mathbb{Q}_{32}}(f)]$
is not zero, and we have proved case (A) of Theorem \ref{th:Main1}.

\subsubsection{Case $n=10$ and $\mathbf{\protect\alpha }=(1,2,2,1,2,2)$}

\noindent Since the proof of (B) follows the steps of the previous case, we
will just outline the computational parts which differ.

\noindent Let $f:S^{3}\rightarrow W_{10}$ be given by $f(t)=(-\frac{%
21809669044}{3234846615},\frac{23}{29},\frac{19}{23},\frac{17}{19},\frac{13}{%
17},\frac{11}{13},\frac{7}{11},\frac{5}{7},\frac{3}{5},\frac{2}{3})$.

\noindent The arrangement $\mathcal{A}$ is now a minimal $\mathbb{Q}_{40}$%
-arrangement containing the subspace $L$ defined by%
\begin{equation*}
\begin{array}{cccc}
x_{1}+x_{2}+x_{3}+x_{4}+x_{5} & = & x_{2}+x_{3}+x_{4}+x_{5}+x_{6}= &
x_{3}+x_{4}+x_{5}+x_{6}+x_{7} \\
& = & \sum_{i=1}^{8}x_{i}=0\text{.} &
\end{array}%
\end{equation*}%
The arrangement $\mathcal{A}$ has five maximal elements $L$, $\epsilon L$, $%
\epsilon ^{2}L$, $\epsilon ^{3}L$ and $\epsilon ^{4}L$. This follows from
the set equality $L=\epsilon ^{5}L$. Let us determine the intersection of
the image under $f$ of the maximal cell%
\begin{equation*}
e=\left( [t,\epsilon t]\cup \lbrack \epsilon t,\epsilon ^{2}t]\cup \lbrack
\epsilon ^{2}t,\epsilon ^{3}t]\cup ~...~\cup \lbrack \epsilon ^{8}t,\epsilon
^{9}t]\right) \ast \lbrack jt,\epsilon jt]
\end{equation*}%
with the test space $\bigcup \mathcal{A}=L\cup \epsilon L\cup \epsilon
^{2}L\cup \epsilon ^{3}L\cup \epsilon ^{4}L$.\ The results of $10\times 5=50$
intersections $f\left( [\epsilon ^{i}t,\epsilon ^{i+1}t]\ast \lbrack
jt,\epsilon jt]\right) \cap \epsilon ^{r}L$ can be summarized in the
following way:%
\begin{equation*}
\begin{array}{ccc}
f\left( [\epsilon ^{3}t,\epsilon ^{4}t]\ast \lbrack jt,\epsilon jt]\right)
\cap \epsilon ^{2}L=\{q_{1}\}, &  & f\left( [\epsilon ^{5}t,\epsilon
^{6}t]\ast \lbrack jt,\epsilon jt]\right) \cap \epsilon ^{3}L=\{q_{2}\}, \\
f\left( [\epsilon ^{5}t,\epsilon ^{6}t]\ast \lbrack jt,\epsilon jt]\right)
\cap \epsilon ^{4}L=\{q_{3}\} &  &
\end{array}%
\end{equation*}%
and consequently $\mathrm{card}\left( f(e)\cap \bigcup \mathcal{A}\right) =3$%
. There are $s_{1},s_{2},s_{3}\in \{1,-1\}$ such that%
\begin{equation}
\mathfrak{o}_{\mathbb{Q}_{40}}(f)(e)=s_{1}\left\Vert q_{1}\right\Vert
+s_{2}\left\Vert q_{2}\right\Vert +s_{3}\left\Vert q_{3}\right\Vert .
\label{ObCocyceCase3}
\end{equation}

\noindent As in the previous case we use the map $\chi ^{\ast }$ (\ref%
{projection}) and prove that $\chi ^{\ast }([\mathfrak{o}_{\mathbb{Q}%
_{40}}(f)])\neq 0$. The set equality $L=\epsilon ^{5}L$ implies that%
\begin{equation}
\mathrm{Hom}\left( \underset{V\in P:\dim V=6}{\bigoplus }\tilde{H}%
_{-1}(\Delta (P_{<V});%
\mathbb{Z}
\mathbf{),}%
\mathbb{Z}
\right) _{\mathbb{Q}_{40}}~\mathbb{\cong ~}%
\mathbb{Z}
_{2}.  \label{iso-1}
\end{equation}%
The element $\epsilon ^{5}$ acts on $W_{10}$ by changing its orientation. On
the orthogonal complement $L^{\perp }$ of $L$ the operator $\epsilon ^{5}$,
for the basis $%
\{e_{1}+...+e_{5},~e_{2}+...+e_{6},~e_{3}+...+e_{7},~e_{1}+...+e_{10}\}$ of $%
L^{\perp }$, has the matrix%
\begin{equation*}
\Xi =\left(
\begin{array}{cccc}
{\small -1} & {\small 0} & {\small 0} & {\small 1} \\
{\small 0} & {\small -1} & {\small 0} & {\small 1} \\
{\small 0} & {\small 0} & {\small -1} & {\small 1} \\
{\small 0} & {\small 0} & {\small 0} & {\small 1}%
\end{array}%
\right) \text{.}
\end{equation*}%
Since $\det \Xi =-1$, the element $\epsilon ^{5}$ changes the orientation of
$L^{\perp }$ and consequently does not change the orientation on $L$. Let $%
l\in H_{6}(\bigcup \widehat{\mathcal{A}};%
\mathbb{Z}
)$ be the generator associated to the subspace $L$, then the sum $\underset{%
V\in P:\dim V=6}{\bigoplus }\tilde{H}_{-1}(\Delta (P_{<V});\mathbb{%
\mathbb{Z}
}\mathbf{)}$ is generated as a free abelian group by elements $l$, $\epsilon
l$, $\epsilon ^{2}l$, $\epsilon ^{3}l$, $\epsilon ^{4}l$. The set equality $%
L=\epsilon ^{5}L$ implies the homology equality $\epsilon ^{5}l=l$. Let $\xi
\in $ $\mathrm{Hom}\left( \underset{V\in P:\dim V=4}{\bigoplus }\tilde{H}%
_{-1}(\Delta (P_{<V});\mathbb{%
\mathbb{Z}
}\mathbf{)},\mathbb{%
\mathbb{Z}
}\right) $ be defined by $\xi (l)=1$ and $\xi (\epsilon ^{i}l)=0$ for $i\in
\{1,2,3,4\}$. Then
\begin{equation*}
(\epsilon ^{5}\ast \xi )(l)=\det (\epsilon ^{5})\xi (\epsilon ^{-5}l)=-\xi
(\epsilon ^{5}l)=-\xi (l)=-1
\end{equation*}%
and
\begin{equation*}
(\epsilon ^{5}\ast \xi )(\epsilon ^{i}l)=\det (\epsilon ^{5})\xi (\epsilon
^{-5+i}l)=-\xi (\epsilon ^{-5+i}l)=0
\end{equation*}%
for $i\in \{1,2,3,4\}$. Thus, there is a relation
\begin{equation*}
\xi +\epsilon ^{5}\ast \xi =0
\end{equation*}%
which implies the isomorphism (\ref{iso-1}) we indicated. As before $\chi
^{\ast }([\mathfrak{o}_{\mathbb{Q}_{40}}(f)])$ is a generator of the group $%
\mathbb{Z}
_{2}$. Therefore case (B) of Theorem \ref{th:Main1} is proved.

\noindent We state the following conjecture.

\begin{conjecture}
All symmetric six-tuples $(a,b,c,a,b,c)$ are solutions of problem \ref%
{Problem1}.
\end{conjecture}

\section{\label{Result2}The $(a,b,a)$ class of $3$-fan $2$-measures
partitions}

\noindent The problem discussed in this chapter will be managed by an
extension of obstruction theory methods used in the paper \cite{Bl}. Our
result, although obtained along the lines of the same method, differs from
the result of the paper \cite{Bl} at the most critical phase of the proof.
The search for the target extension space and the identification of the
obstruction element is only remotely similar.

\noindent The appendix contains the ancillary material referred to in this
section.

\bigskip

\noindent Let $\mu _{1},\mu _{2},\ldots ,\mu _{m}$ be proper Borel
probability measures on $S^{2}$. Let $(\alpha _{1},\alpha _{2},\ldots
,\alpha _{k})\in \left( \mathbb{R}_{>0}\right) ^{k}$ be a vector where $%
\alpha _{1}+\alpha _{2}+\ldots +\alpha _{k}=1$. The general problem stated
in \cite{BaMa2001} is:

\begin{problem}
Determine all triples $(m,k,\alpha )\in \mathbb{N\times }\mathbb{N\times R}%
^{k}$ such that for any collection of $m$ measures $\{\mu _{1},\mu
_{2},\ldots ,\mu _{m}\}$, there exists a $k$-fan $(x;l_{1},\ldots ,l_{k})$
with the property%
\begin{equation*}
(\forall i=1,\ldots ,k)\,(\forall j=1,\ldots ,m)\,\mu _{j}(\sigma
_{i})=\alpha _{i}.
\end{equation*}%
Such a $k$-fan $(x;l_{1},\ldots ,l_{k})$ is called an $\alpha $\textit{%
-partition} for the collection of measures $\{\mu _{1},\mu _{2},\ldots ,\mu
_{m}\}$.
\end{problem}

\noindent In this chapter we prove the existence of a class of solutions of
the above problem in the case $m=2$ / $k=3$.

\begin{theorem}
\label{th:main2}Let us choose $\alpha =(a,b,a)\in
\mathbb{R}
_{>0}^{3}$ such that $2a+b=1$. Then any two \textit{proper} Borel
probability measures $\mu $ and $\nu $ on the sphere $S^{2}$ admit an $%
\alpha $-partition by a $3$-fan $\mathfrak{p}=(x;l_{1},l_{2},l_{3})$.
\end{theorem}

\subsection{The configuration space / test map scheme}

\noindent The CS / TM setting is the same as in the previous chapter except
for the test space. Therefore we give only an outline. For details one can
also consult \cite{BaMa2001}, \cite{Bl} and \cite{Bl-Vr-Ziv}.

\noindent \textbf{The configuration space.} Let $\mu $ and $\nu $ be two
proper Borel probability measures on $S^{2}$. The configuration space $%
X_{\mu }$ associated with the measure $\mu $ is
\begin{equation*}
X_{\mu }=\{(x;l_{1},\ldots ,l_{n})\in F_{n}\mid (\forall i=1,\ldots ,n)\,\mu
(\sigma _{i})=\tfrac{{\small 1}}{{\small n}}\}.
\end{equation*}%
As before, $X_{\mu }$ is a Stiefel manifold $V_{2}(%
\mathbb{R}
^{3})$ of all orthonormal $2$-frames in $%
\mathbb{R}
^{3}$.

\noindent \textbf{Test maps.} Let $\alpha $-vectors have the form $\alpha =(%
\tfrac{a_{1}}{n},\tfrac{a_{2}}{n},\tfrac{a_{3}}{n})\in \frac{1}{n}\,\mathbb{N%
}^{3}$, where $a_{1}+a_{2}+a_{3}=n$. The test map for this fan problem is
defined by%
\begin{equation*}
\begin{array}{l}
\Phi _{\nu }:X_{\mu }\rightarrow W_{n}=\{x\in
\mathbb{R}
^{n}\mid x_{1}+\ldots +x_{n}=0\} \\
\Phi _{\nu }(\mathfrak{p})=\left( \nu (\sigma _{1})-\tfrac{{\small 1}}{%
{\small n}},\ldots ,\nu (\sigma _{n})-\tfrac{{\small 1}}{{\small n}}\right) .%
\end{array}%
\end{equation*}%
\noindent \textbf{The action.} The dihedral group $\mathbb{D}_{2n}=\langle
j,\varepsilon \,|\,\varepsilon ^{n}=j^{2}=1,\,\varepsilon j=j\varepsilon
^{n-1}\,\rangle $ acts both on the possible solution space $X_{\mu }$ and
the linear subspace $W_{n}\subseteq
\mathbb{R}
^{n}$, as in the previous chapter. Observe that $X_{\mu }$ is $\mathbb{D}%
_{2n}$-homeomorphic to the manifold $V_{2}(%
\mathbb{R}
^{3})$, where $V_{2}(%
\mathbb{R}
^{3})$ is a $\mathbb{D}_{2n}$-space given by
\begin{equation*}
\varepsilon (x,y)=(x,R_{x}({\tfrac{2\pi }{n}})(y))\text{, }j(x,y)=(-x,y)%
\text{,}
\end{equation*}%
and $R_{x}(\theta ):%
\mathbb{R}
^{3}\rightarrow
\mathbb{R}
^{3}$ is the rotation around the axes determined by $x$ through the angle $%
\theta $.

\noindent \textbf{The test space. }Given $\alpha =(\frac{a_{1}}{n},\frac{%
a_{2}}{n},\frac{a_{3}}{n})\in \frac{1}{n}\,\mathbb{N}^{3}$, let $L(\alpha
)\subset W_{n}$ be defined by the equations%
\begin{equation}
x_{1}+\ldots +x_{a_{1}}=x_{a_{1}+1}+\ldots
+x_{a_{1}+a_{2}}=x_{a_{1}+a_{2}+1}+\ldots +x_{n}=0\text{.}  \label{L(alfa)}
\end{equation}%
The test space for this problem is the union $\bigcup \mathcal{A}(\alpha
)\subset W_{n}$ of the smallest $\mathbb{D}_{2n}$-invariant linear subspace
arrangement $\mathcal{A}(\alpha )$, containing $L(\alpha )$.

\noindent Since the test map $\Phi _{\nu }$ is $\mathbb{D}_{2n}$%
-equivariant, the following proposition holds.

\begin{proposition}
\label{prop:FansVezaProblem-Ekvivarijantan}Let $\alpha =(\frac{a_{1}}{n},%
\frac{a_{2}}{n},\frac{a_{3}}{n})\in \frac{1}{n}\,\mathbb{N}^{3}$ be a vector
such that $a_{1}+a_{2}+a_{3}=n$. If there is no $\mathbb{D}_{2n}$%
-equivariant map%
\begin{equation*}
\Phi :V_{2}(%
\mathbb{R}
^{3})\rightarrow W_{n}\setminus \bigcup \mathcal{A}(\alpha )
\end{equation*}%
then for any two measures $\mu $ and $\nu $ on $S^{2}$, there exists an $%
\alpha $-partition $(x;l_{1},l_{2},l_{3})$ of measures $\mu $ and $\nu $.
\end{proposition}

\noindent As we have seen in Proposition \ref{prop:PrelazakNaNovuGrupu} and
in \cite{Bl} and \cite{Bl-Vr-Ziv} the Stiefel manifold $V_{2}(%
\mathbb{R}
^{3})$ can be substituted with the sphere $S^{3}$. Thus we prove that there
is no $\mathbb{Q}_{4n}$-map%
\begin{equation*}
S^{3}\rightarrow W_{n}\setminus \bigcup \mathcal{A}(\alpha ).
\end{equation*}

\subsection{Proof of Theorem \protect\ref{th:main2}}

\noindent Theorem \ref{th:main2} will be proved by showing that there is no $%
\mathbb{Q}_{4n}$-map $S^{3}\rightarrow W_{n}\setminus \bigcup \mathcal{A}%
(a,b,a)$. The arrangement $\mathcal{A}(a,b,a)$ is a minimal $\mathbb{D}_{2n}$%
-invariant arrangement containing the subspace $L\subset W_{n}$ given by
\begin{equation*}
x_{1}+\ldots +x_{a}=x_{a+1}+\ldots +x_{a+b}=x_{a+b+1}+\ldots +x_{n}=0\text{.}
\end{equation*}%
The codimension of the arrangement $\mathcal{A}(a,b,a)$ inside $W_{n}$ is
two, and so the complement is connected. The fundamental group of the
complement is far from trivial, so the obstruction theory can not be applied
directly. We use the idea of a target extension scheme introduced in \cite%
{Bl}. Let us denote by $M$ the complement $W_{n}\setminus \bigcup \mathcal{A}%
(a,b,a)$.

\subsubsection{The target extension scheme.}

The basic idea of the target extension scheme is to find a $\mathbb{Q}_{4n}$%
-space $N$\ which contains the target space $M$, and to prove that there is
no $\mathbb{Q}_{4n}$-map\emph{\ }$S^{3}\rightarrow N$. This would imply that
there is no $\mathbb{Q}_{4n}$-map\textit{\ }$S^{3}\rightarrow M$. In our
case the target space $M$ is the complement of the arrangement, and so we
can refine the idea as follows:

\begin{itemize}
\item \emph{Increase the codimension.}\textit{\ }Let $H$ be an arbitrary
hyperplane in $W_{n}$ and let $\mathcal{B}$ be the minimal $\mathbb{Q}_{4n}$%
-invariant arrangement containing the subspace $L\cap H$. The inclusion $%
\bigcup \mathcal{B}\subseteq \bigcup \mathcal{A}(a,b,a)$ implies that
\begin{equation}
W_{n}\backslash \bigcup \mathcal{B}\supseteq W_{n}\backslash \bigcup
\mathcal{A}(a,b,a).
\end{equation}%
The dimension of maximal elements of the arrangement $\mathcal{B}$ is $n-4$.
Let us denote by $N$ the "new" complement $W_{n}\backslash \bigcup \mathcal{B%
}$.

\item \emph{Apply the general position map scheme.} The codimension of the
arrangement $\mathcal{B}$ inside $W_{n}$ is three and so the complement $%
N=W_{n}\backslash \bigcup \mathcal{B}$ is 1-connected. Therefore the
question of the existence of a $\mathbb{Q}_{4n}$-map\textit{\ }$%
S^{3}\rightarrow N$ depends only on the primary obstruction. We can now use
the general position map scheme.
\end{itemize}

\noindent Unfortunately, the target extension scheme only provides hope that
we can apply, once more, the mechanisms we already developed. The main
problem is just shifted to the question of \emph{how to find a hyperplane} $H
$ in a way that there is no $\mathbb{Q}_{4n}$-map $S^{3}\rightarrow N$.

\subsubsection{In the pursuit of the hyperplane $H$.}

There are two rough heuristics we can use. First, if we define a $\mathbb{Q}%
_{4n}$-map\textit{\ }$f:S^{3}\rightarrow W_{n}$, we can introduce $H$ in
such a way that

(A) $f$ becomes a map in general position for the new arrangement $\mathcal{B%
}$, and

(B) the cardinality of the set $f(S^{3})\cap \bigcup \mathcal{B}$ is as
small as can be achieved.

\noindent The second requirement for introduction of the hyperplane is that
the group of coinvariants $H_{2}(W_{n}\backslash \bigcup \mathcal{B};\mathbb{%
Z})_{\mathbb{Q}_{4n}}$ has torsion. Nevertheless, the choice of the
hyperplane $H$ is strongly connected with the properties of the action of
the group on the arrangement $\mathcal{A}(a,b,a)$. In this situation we are
going to exploit the following symmetry $jL=L$. Observe that element $j\in
\mathbb{Q}_{4n}$ \emph{changes} the orientation of the orthogonal complement
$L^{\bot }$.

\subsubsection{The map in general position.}

Let us define a $\mathbb{Q}_{4n}$-map\textit{\ }$f:S^{3}\rightarrow W_{n}$
and introduce $H$ so that $f$ becomes a map in general position. Let $S^{3}$
be the simplicial complex $P_{2n}\ast P_{2n}$, where $P_{2n}$ is the regular
$2n$-gon. Let $t$ be a fixed vertex in one of the copies of $P_{2n}$. Let $%
u_{i}=e_{i}-\tfrac{1}{n}\sum\limits_{j=1}^{n}e_{j}$, $i\in \{1,..,n\}$,
where $e_{1},..,e_{n}$ are elements of the standard basis of $\mathbb{R}^{n}$%
. We define $f:S^{3}\rightarrow W_{n}$ on the vertex $h(t)=u_{1}$ and then
extend it equivariantly. This implies
\begin{equation*}
\begin{array}{lll}
h(\epsilon ^{i}t)=\epsilon ^{i}h(t)=u_{\left( i+1\right) ~\mathrm{mod}~n} &
& h(\epsilon ^{i}jt)=\epsilon ^{i}jh(t)=u_{i~\mathrm{mod}~n}.%
\end{array}%
\end{equation*}%
In the future all the indexes in $W_{n}$ will be calculated $\mathrm{mod}~n$.

\begin{notation}
Every 3-simplex of the sphere $P_{2n}\ast P_{2n}$, determined by four
vertices, has form
\begin{equation}
\lbrack \epsilon ^{p}t,\epsilon ^{p+1}t;\epsilon ^{q}jt,\epsilon ^{q+1}jt]
\end{equation}%
and will be denoted by $\sigma _{p,q}$.
\end{notation}

\noindent The list of all simplices of the form $%
[u_{i},u_{i+1};u_{j},u_{j+1}]$ which intersect the subspace $L$ is given in
the following table

\begin{equation*}
\begin{tabular}{|l|l|l|}
\hline
{\small \ \ Simplex} & {\small Bar. coor. of the intersection } & {\small %
for }${\small r\in }$ \\ \hline
{\small 1. \ }$[u_{a},u_{a+1};u_{a+b},u_{a+b+1}]$ & $\frac{1}{n}{\small %
(a,\beta ,\gamma ,a);~\beta +\gamma =b}$ &  \\ \hline
{\small 2. \ }$[u_{a},u_{a+1};u_{r},u_{r+1}]$ & $\frac{1}{n}{\small %
(a,b,\gamma ,\delta );~\gamma +\delta =a}$ & ${\small [a+b+1,n-3]}$ \\ \hline
{\small 3. \ }$[u_{a},u_{a+1};u_{n-2},u_{n-1}]$ & $\frac{1}{n}{\small %
(a,b,\gamma ,\delta );~\gamma +\delta =a}$ &  \\ \hline
{\small 4. \ }$[u_{a},u_{a+1};u_{n-1},u_{n}]$ & $\frac{1}{n}{\small %
(a,b,\gamma ,\delta );~\gamma +\delta =a}$ &  \\ \hline
{\small 5. \ }$[u_{r},u_{r+1};u_{a+b},u_{a+b+1}]$ & $\frac{1}{n}{\small %
(\alpha ,\beta ,b,a);~\alpha +\beta =a}$ & ${\small [2,a-2]}$ \\ \hline
{\small 6. \ }$[u_{1},u_{2};u_{a+b},u_{a+b+1}]$ & $\frac{1}{n}{\small %
(\alpha ,\beta ,b,a);~\alpha +\beta =a}$ &  \\ \hline
{\small 7. \ }$[u_{n},u_{1};u_{a},u_{a+1}]$ & $\frac{1}{n}{\small (\alpha
,\beta ,b,a);~\alpha +\beta =a}$ &  \\ \hline
{\small 8. \ }$[u_{n},u_{1};u_{a+b},u_{a+b+1}]$ & $\frac{1}{n}{\small %
(a,b,\gamma ,\delta );~\gamma +\delta =a}$ &  \\ \hline
{\small 9. \ }$[u_{n},u_{1};u_{r},u_{r+1}]$ & $\frac{1}{n}{\small (a,\beta
,\gamma ,a);~\beta +\gamma =b}$ & ${\small [a+1,a+b-2]}$ \\ \hline
{\small 10. }$[u_{n},u_{1};u_{a+b-1},u_{a+b}]$ & $\frac{1}{n}{\small %
(a,\beta ,\gamma ,a);~\beta +\gamma =b}$ &  \\ \hline
\end{tabular}%
\end{equation*}%
If we introduce the hyperplane $H$ by
\begin{equation*}
H=\{\mathbf{x}\in
\mathbb{R}
^{n}|~x_{a+b}-x_{a+1}+x_{1}-x_{n-1}=0\},
\end{equation*}%
the list of simplices now intersecting $L_{1}=L\cap H$ reduces to

(A) two simplices, for $a>b$,
\begin{equation*}
\rho _{1}=[u_{a},u_{a+1};u_{a+b},u_{a+b+1}]\text{ and }\rho
_{2}=[u_{n},u_{1};u_{a},u_{a+1}],
\end{equation*}%
with barycentric coordinates of intersection points $\tfrac{1}{n}{\small (a,}%
\tfrac{{\small b}}{2}{\small ,\tfrac{b}{2},a)}$ and $\tfrac{1}{n}{\small %
(b,a-b,b,a);}$

(B) one simplex, for $b>a$, $\rho _{1}=[u_{a},u_{a+1};u_{a+b},u_{a+b+1}]$,
with barycentric coordinates of intersection points $\tfrac{1}{n}{\small (a,}%
\tfrac{{\small b}}{2}{\small ,\tfrac{b}{2},a)}$.

\noindent In order to simplify the exposition let us assume that $b>a$. Case
(A) is actually a consequence of Case (B). Therefore,%
\begin{equation*}
f(S^{3})\cap L_{1}=\{\mathbf{y}=\tfrac{a}{n}u_{a}+\tfrac{b}{2n}u_{a+1}+%
\tfrac{b}{2n}u_{a+b}+\tfrac{a}{n}u_{a+b+1}\}.
\end{equation*}

\noindent We have made extensive use of the first heuristic. In is not easy
to see that we also used the second heuristic. Let us just say that the
following property will be a significant ingredient in the interpretation of
the obstruction element.%
\begin{equation}
\text{Element }j\text{ does NOT change the orientation of the subspace }%
(L_{1}\cap jL_{1})^{\bot }.  \label{j(L1)}
\end{equation}

\subsubsection{The obstruction cocycle.}

\noindent There are 8 simplices in $P_{2n}\ast P_{2n}$ which are in the
inverse image under $f$ of the simplex $[u_{a},u_{a+1};u_{a+b},u_{a+b+1}]$.
These simplices are%
\begin{equation*}
\begin{array}{llll}
\sigma _{a-1,a+b}; & \sigma _{n+a-1,a+b}; & \sigma _{a-1,n+a+b}; & \sigma
_{n+a-1,n+a+b}; \\
\sigma _{a+b-1,a}; & \sigma _{n+a+b-1,a}; & \sigma _{a+b-1,n+a}; & \sigma
_{n+a+b-1,n+a}.%
\end{array}%
\end{equation*}%
All eight simplices are in the orbit of the following two simplices from the
maximal cell $e$ (consult appendix \ref{Sec:Simplicial}):%
\begin{equation*}
\Theta _{1}=\sigma _{2a-1,0}\text{ and }\Theta _{2}=\sigma _{b-1,0}\text{.}
\end{equation*}%
The points $\xi _{1}\in \Theta _{1}$ and $\xi _{2}\in \Theta _{2}$ ,with
barycentric coordinates $(a,\frac{b}{2},\frac{b}{2},a)$ and $(\frac{b}{2}%
,a,a,\frac{b}{2})$, are mapped inside the arrangement $\mathcal{B}$,
precisely%
\begin{equation*}
f(\xi _{1})\in \epsilon ^{a}(L_{1}\cap jL_{1})\text{ and }f(\xi _{2})\in
\epsilon ^{-a}(L_{1}\cap jL_{1})\text{.}
\end{equation*}%
Therefore, formula (\ref{ObCocycle}) implies that the obstruction cocycle is
\begin{equation*}
\begin{array}{lll}
\mathfrak{o}_{\mathbb{Q}_{4n}}(f)(e) & =\left\Vert f(\xi _{1})\right\Vert
+\left\Vert f(\xi _{2})\right\Vert =\det (\epsilon ^{a})\epsilon
^{a}\left\Vert \mathbf{y}\right\Vert +\det (\epsilon ^{-a})\epsilon
^{-a}\left\Vert \mathbf{y}\right\Vert  &  \\
& =\det (\epsilon ^{-a})\left( \det (\epsilon ^{2a})\epsilon ^{a}\left\Vert
\mathbf{y}\right\Vert +\epsilon ^{-a}\left\Vert \mathbf{y}\right\Vert
\right)  &
\end{array}%
\end{equation*}%
Since we are interested in the cohomology / coinvariant class of the
obstruction cocycle (\ref{CoinvariantIsomorphism}), instead of $\mathfrak{o}%
_{\mathbb{Q}_{4n}}(f)(e)$ we can look at the element
\begin{equation}
\mathfrak{o}_{\mathbb{Q}_{4n}}^{\prime }=2\det (\epsilon ^{-a})\left\Vert
\mathbf{y}\right\Vert \text{.}  \label{obCocycleFans}
\end{equation}

\begin{remark}
In contrast to the previous chapter, the point class $\left\Vert \mathbf{y}%
\right\Vert $ \textbf{depends} on the placement of the simplex $\rho
_{1}=[u_{a},u_{a+1};u_{a+b},u_{a+b+1}]$ with respect to the arrangement $%
\mathcal{B}$. In terms of paper \cite{Bl}, $\left\Vert \mathbf{y}\right\Vert
$ is the broken point class.
\end{remark}

\subsubsection{The obstruction element.}

\noindent It remains to prove that $2\left\Vert \mathbf{y}\right\Vert $ does
not vanish when we pass to coinvariants. The brief proof of this fact would
be: There is a symmetry%
\begin{equation*}
j[u_{a},u_{a+1};u_{a+b},u_{a+b+1}]=[u_{a},u_{a+1};u_{a+b},u_{a+b+1}]
\end{equation*}%
which does not change the orientation of the simplex and therefore gives the
relation $j\left\Vert \mathbf{y}\right\Vert =\left\Vert \mathbf{y}%
\right\Vert $ and NOT $\left\Vert \mathbf{y}\right\Vert =-j\left\Vert
\mathbf{y}\right\Vert $. This is the only relation which could provide in
coinvariants that $2(class$ $\left\Vert \mathbf{y}\right\Vert )=0$. Thus $%
2\left\Vert \mathbf{y}\right\Vert $ does not vanish in coinvariants and the
obstruction element is NOT zero. We proved Theorem \ref{th:main2} for all
rational triples $(a,b,a)$.

\bigskip

\noindent Just to be safe, let us analyze the point class $\left\Vert
\mathbf{y}\right\Vert $, its placement with the respect to the arrangement $%
\mathcal{B}$ and the indicated $j$ symmetry. Observe that we are not forced
to determine the complete group of coinvariants $H_{2}\left( W_{n}\backslash
\bigcup \mathcal{B};\mathbb{Z}\right) _{\mathbb{Q}_{4n}}$. We only have to
prove that $2(class$ $\left\Vert \mathbf{y}\right\Vert )\neq 0$.

\bigskip

\begin{figure}[tbh]
\centering \includegraphics[scale=0.70]{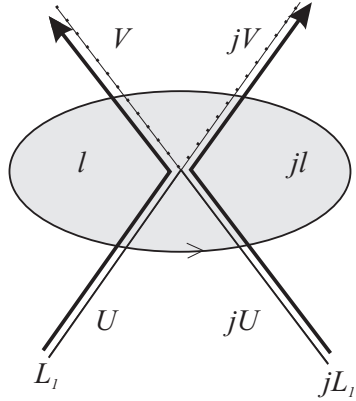}
\caption{{\protect\small The position of a point class relative to the
arrangement}}
\label{fig:Fig2}
\end{figure}

\noindent Translating the simplex $[u_{a},u_{a+1};u_{a+b},u_{a+b+1}]$ by a
small generic vector and then intersecting it with the arrangement, we find
that the boundary (of the shaded region) links with the elements of the
arrangement as in the Figure \ref{fig:Fig2}.

\noindent Figure \ref{fig:Fig2} \ indicates two set relations $L_{1}=U\cup jV
$ and $jL_{1}=jU\cup V$. Since element $j$ does not change the orientation
of the subspace $(L_{1}\cap jL_{1})^{\bot }$, the geometric generators of
the group $H_{n-4}(\bigcup \widehat{\mathcal{B}};\mathbb{%
\mathbb{Z}
})$ (the dual of $H_{2}(W_{n}\backslash \bigcup \mathcal{B};\mathbb{Z})$ )
can only be given by sets%
\begin{equation}
L_{1},\epsilon L_{1},..,\epsilon ^{n-1}L_{1},jL_{1},\epsilon
jL_{1},..,\epsilon ^{n-1}jL_{1}\text{ and }U\cup V,\epsilon \left( U\cup
V\right) ,..,\epsilon ^{n-1}\left( U\cup V\right) \text{.}
\label{homologyBasis}
\end{equation}%
This fact is a consequence of the following observations

\begin{itemize}
\item the orientation on the arrangement should be consistent with the group
action,

\item the boundaries of the two pieces which are glued together have
opposite orientations,

\item the element $j\left( U\cup V\right) $ is not a generator because it is
a linear combination of already introduced generators.
\end{itemize}

\noindent The isomorphism composition of the Poincar\'{e} duality
isomorphism and the Universal coefficient isomorphism for the arrangement $%
\mathcal{B}$,
\begin{equation*}
\varphi :H_{2}\left( W_{n}\setminus \bigcup \mathcal{B}\right) \rightarrow
\mathrm{Hom}\left( H_{n-4}\left( \bigcup \widehat{\mathcal{B}},\mathbb{%
\mathbb{Z}
}\right) ;\mathbb{%
\mathbb{Z}
}\right) ,
\end{equation*}%
is the computation of the linking number with the elements of the basis (\ref%
{homologyBasis}). Thus,%
\begin{equation*}
\begin{tabular}{ccc}
$\varphi (\left\Vert \mathbf{y}\right\Vert )(l)=\pm 1,$ & $\varphi
(\left\Vert \mathbf{y}\right\Vert )(jl)=\mp 1,$ & $\varphi (\left\Vert
\mathbf{y}\right\Vert )(t)=1,$ \\
$\varphi (\left\Vert \mathbf{y}\right\Vert )(u)=0,$ & \multicolumn{2}{c}{for
all other basis elements,}%
\end{tabular}%
\end{equation*}%
where $l,t\in $ $H_{n-5}(\bigcup \widehat{\mathcal{B}},\mathbb{%
\mathbb{Z}
)}$ are elements of a basis determined by sets $L_{1}$ and $U\cup V$. Since $%
\varphi (\left\Vert \mathbf{y}\right\Vert )(t)=1$, the class $2[\varphi
(\left\Vert \mathbf{y}\right\Vert )]$ does not vanish. The obstruction
element is NOT zero and we have proved Theorem \ref{th:main2} for all
rational triples $(a,b,a)$.

\subsubsection{The closing argument.}

\noindent Let $\emph{S}\subseteq \mathbb{R}_{>0}^{3}$ be the space of all
triples $(a,b,c)$, $a+b+c=1$ such that for $\alpha =(a,b,c)$ there exists an
$\alpha $-partition of measures $\mu $ and $\nu $ by a $3$-fan. Since we
assumed that our measures $\mu $ and $\nu $ are two proper Borel probability
measures, the space $\emph{S}$ is a closed subset of $\mathbb{R}_{>0}^{3}$.
In this notation, so far we have proved the inclusion
\begin{equation*}
\{(\tfrac{a}{n},\tfrac{b}{n},\tfrac{a}{n})\in \mathbb{Q}%
_{>0}^{3}~|~2a+b=n,~a,b\in \mathbb{Z}\}\subseteq \emph{S.}
\end{equation*}%
Therefore
\begin{equation*}
\{(a,b,c)\in \mathbb{R}_{>0}^{3}~|~a+b+c=1\}=\mathrm{cl}\{(\tfrac{a}{n},%
\tfrac{a+b}{n},\tfrac{b}{n})~|~2a+2b=n,~a,b\in \mathbb{Z}\}\subseteq \emph{S}%
\text{.}
\end{equation*}%
where $\mathrm{cl}$ denotes the closure of the subset in $\mathbb{R}^{3}$.
The theorem is proved.

\section{Appendix: Geometry of $\mathbb{Q}_{4n}$}

\subsection{\label{Sec:Q4n}Generalized quaternion group.}

The sphere $S^{3}=S(\mathbb{H})=Sp(1)$ can be seen as the group of all unit
quaternions. Consider a root of unity $\epsilon =\epsilon _{2n}=\cos \frac{%
\pi }{n}+i\sin \frac{\pi }{n}\in S(\mathbb{H})$. The group generated by $%
\epsilon $ is a subgroup of $S(\mathbb{H})$ of order $2n$. The generalized
quaternion group\textit{, }\cite[p. 253]{CaEi}\textit{, }is a subgroup of
order $4n$ generated by $\epsilon $ and $j$. The group $\mathbb{Q}_{4n}$
acts on $S^{3}$ as a subgroup, and on $W_{n}$ via the already defined $%
\mathbb{D}_{2n}$-action by the quotient homomorphism $\mathbb{Q}%
_{4n}\rightarrow \mathbb{Q}_{4n}/\{1,-1\}\cong \mathbb{D}_{2n}$. The $%
\mathbb{Q}_{4n}$-action on $S^{3}$ is \textbf{free}.

\noindent Let $H=$ $\{1,\epsilon ^{n}\}=\{1,-1\}\subset \mathbb{Q}_{4n}$.
Then the quotient group $\mathbb{Q}_{4n}/H$ is isomorphic to the dihedral
group $\mathbb{D}_{2n}$ of order $2n$.\textbf{\ }Also, the group $H$ acts on
$W_{n}$ and on $%
\mathbb{R}
^{n}$ trivially.

\noindent Since the group $\mathbb{Q}_{4n}$ acts on $S^{3}$ as a subgroup
and $S^{3}$ is a connected group, then the $\mathbb{Q}_{4n}$-module
structure $\mathcal{Z}$ on $H_{3}(S^{3},%
\mathbb{Z}
)$ must be trivial.

\subsection{\label{Sec:Simplicial}The natural $\mathbb{Q}_{4n}$ cellular /
simplicial structure on $S^{3}$.}

Let the circle $S^{1}$ be represented by the simplicial complex of the
regular $2n$-gon $P_{2n}$. Then the sphere $S^{3}$, as a simplicial complex,
is the join $P_{2n}^{(1)}\ast P_{2n}^{(2)}$ of two copies of $P_{2n}$.
Denote the vertices of $P_{2n}^{(1)}$ and $P_{2n}^{(2)}$ by $a_{1},..,a_{2n}$
and $b_{1},..,b_{2n}$, respectively. The action of the group $\mathbb{Q}%
_{4n} $ on $S^{3}$ is defined on vertices by%
\begin{equation*}
\epsilon \cdot a_{i}=a_{i~\text{mod}2n+1}\text{, }\epsilon \cdot b_{i}=b_{i~%
\text{mod}2n+1}\text{, }j\cdot a_{1}=b_{1}
\end{equation*}%
and it extends equivariantly to the upper skeletons. The associated chain
complex $\mathfrak{C}=\{\emph{C}_{i}\}$ has the form%
\begin{equation*}
0\longrightarrow \mathbb{%
\mathbb{Z}
}^{4n^{2}}\longrightarrow \mathbb{%
\mathbb{Z}
}^{8n^{2}}\longrightarrow \mathbb{%
\mathbb{Z}
}^{4n^{2}+4n}\longrightarrow \mathbb{%
\mathbb{Z}
}^{4n}\longrightarrow 0.
\end{equation*}

\subsection{An economic $\mathbb{Q}_{4n}$ cellular structure on $S^{3}$%
\textbf{.}}

It comes from the minimal resolution of $\mathbb{%
\mathbb{Z}
}$ by free $\mathbb{Q}_{4n}$-modules described in \cite[p.~253]{CaEi}. The
associated cellular complex has one $\mathbb{Q}_{4n}$ $0$-cell $a$, two $%
\mathbb{Q}_{4n}$ $1$-cells $b$ and $b^{\prime }$, two $\mathbb{Q}_{4n}$ $2$%
-cells $c$ and $c^{\prime }$, and then finally one $\mathbb{Q}_{4n}$ $3$%
-cell $e$. The associated chain complex $\mathfrak{D}=\{\emph{D}_{i}\}$ has
the form
\begin{equation*}
0\rightarrow \mathbb{%
\mathbb{Z}
}[\mathbb{Q}_{4n}]e\overset{\partial }{\rightarrow }\mathbb{%
\mathbb{Z}
}[\mathbb{Q}_{4n}]c\oplus \mathbb{%
\mathbb{Z}
}[\mathbb{Q}_{4n}]c^{\prime }\overset{\partial }{\rightarrow }\mathbb{%
\mathbb{Z}
}[\mathbb{Q}_{4n}]b\oplus \mathbb{%
\mathbb{Z}
}[\mathbb{Q}_{4n}]b^{\prime }\overset{\partial }{\rightarrow }\mathbb{%
\mathbb{Z}
}[\mathbb{Q}_{4n}]e\rightarrow 0
\end{equation*}%
where
\begin{equation*}
\begin{array}{lll}
\partial e=(\epsilon -1)c-(\epsilon j-1)c^{\prime } & \partial
c=(1+...+\epsilon ^{n-1})b-(j+1)b^{\prime } & \partial b=(\epsilon -1)a \\
& \partial c^{\prime }=(\epsilon j+1)b+(\epsilon -1)b^{\prime } & \partial
b^{\prime }=(j-1)a%
\end{array}%
.
\end{equation*}%
Thus, it will be enough to look at the obstruction cocycle $c_{\mathbb{Q}%
_{4n}}(h)$ on the maximal cell $e$, and to prove that its image is or is not
zero, when we pass to the cohomology.

\subsection{The chain map $\mathfrak{D}\rightarrow \mathfrak{C}$\textbf{.}}

There exists a cellular map $\mathfrak{f}:\mathfrak{D}\rightarrow \mathfrak{C%
}$. Since we are interested in $3$-cochains we only need to know the
concrete expression for the cellular map on the top dimensional cell $e$,%
\begin{equation*}
\mathfrak{f}(e)=[a_{1},\epsilon a_{1}]\ast \lbrack b_{1},\epsilon
b_{1}]+[\epsilon a_{1},\epsilon ^{2}a_{1}]\ast \lbrack b_{1},\epsilon
b_{1}]+...+[\epsilon ^{n-1}a_{1},\epsilon ^{n}a_{1}]\ast \lbrack
b_{1},\epsilon b_{1}]\text{,}
\end{equation*}%
where the simplices on the right hand side are appropriately oriented.

\end{document}